\documentclass{amsart}
\usepackage{amsmath, amsthm, amsfonts, hyperref, graphicx, ifpdf, mathrsfs}
\usepackage{amssymb}
\usepackage{faktor}
\usepackage{bm}
\usepackage{enumerate}
\usepackage[utf8]{inputenc}
\usepackage{mathtools}
\usepackage{esint}
\usepackage{tikz-cd} 
\hypersetup{
   unicode=false,          % non-Latin characters in AcrobatÃ∞Â¢Ã¢â•¢Â¢Ã—Ë˚s bookmarks
   pdftoolbar=true,        % show AcrobatÃ∞Â¢Ã¢â•¢Â¢Ã—Ë˚s toolbar?
   pdfmenubar=true,        % show AcrobatÃ∞Â¢Ã¢â•¢Â¢Ã—Ë˚s menu?
   pdffitwindow=false,     % window fit to page when opened
   pdfstartview={FitH},    % fits the width of the page to the window
   pdftitle={My title},    % title
   pdfauthor={Author},     % author
   pdfsubject={Subject},   % subject of the document
   pdfcreator={Creator},   % creator of the document
   pdfproducer={Producer}, % producer of the document
   pdfkeywords={keywords}, % list of keywords
   pdfnewwindow=true,      % links in new window
   colorlinks=true,       % false: boxed links; true: colored links
   linkcolor=blue,          % color of internal links
   citecolor=blue,        % color of links to bibliography
   filecolor=magenta,      % color of file links
   urlcolor=MidnightBlue          % color of external links
}

\theoremstyle{plain}
\newtheorem{theorem}{Theorem}[section]
\newtheorem{lem}[theorem]{Lemma}
\newtheorem{pro}[theorem]{Proposition}
\newtheorem{cor}[theorem]{Corollary}
\newtheorem{corx}{Corollary}
\newtheorem{thx}{Theorem}

\theoremstyle{definition}
\newtheorem{definition}[theorem]{Definition}

\theoremstyle{remark}
\newtheorem{rem}[theorem]{Remark}

\numberwithin{equation}{section}

\newcommand{\coc}{cohomology class}
\newcommand{\cp}{critical point}

\newcommand{\cs}{conformal structure}
\newcommand{\cmli}{conformal minimal Lagrangian immersion}
\newcommand{\cmc}{constant mean curvature}

\newcommand{\df}{the Donaldson functional}

\newcommand{\fg}{fundamental group}

\newcommand{\hem}{Hermitian metric}
\newcommand{\hcd}{holomorphic cubic differential}
\newcommand{\hkd}{holomorphic $k$-differential}
\newcommand{\hqd}{holomorphic quadratic differential}

\newcommand{\hym}{hyperbolic metric}
\newcommand{\htm}{hyperbolic three-manifold}

\newcommand{\kf}{K\"{a}hler form}

\newcommand{\mc}{mean curvature}

\newcommand{\mi}{minimal immersion}

\newcommand{\MS}{moduli space}

\newcommand{\mli}{minimal Lagrangian immersion}

\newcommand{\msq}{minimizing sequence}
\newcommand{\mpa}{mountain-pass}

\newcommand{\PS}{Palais-Smale}
\newcommand{\PSc}{Palais-Smale condition}
\newcommand{\qf}{quasi-Fuchsian}

\newcommand{\RS}{Riemann surface}

\newcommand{\sff}{second fundamental form}

\newcommand{\TS}{Teichm\"{u}ller space}

\newcommand{\tm}{three-manifold}

\newcommand{\ubd}{uniformly bounded}
\newcommand{\upto}{along a subsequence}

\newcommand{\wrt}{with respect to}

\newcommand{\be}{\begin{equation}}
\newcommand{\ene}{\end{equation}}
\newcommand{\br}{\begin{rem}}
\newcommand{\er}{\end{rem}}
\newcommand{\bl}{\begin{lem}}
\newcommand{\bcor}{\begin{cor}}
\newcommand{\ecor}{\end{cor}}
\newcommand{\el}{\end{lem}}
\newcommand{\bd}{\begin{Def}}
\newcommand{\ed}{\end{Def}}
\newcommand{\ben}{\begin{enumerate}}
\newcommand{\een}{\end{enumerate}}
\newcommand{\bp}{\begin{proof}}
\newcommand{\ep}{\end{proof}}
\newcommand{\bpo}{\begin{pro}}
\newcommand{\epo}{\end{pro}}
\newcommand{\beq}{\begin{equation*}}
\newcommand{\eeq}{\end{equation*}}
\newcommand{\bear}{\begin{eqnarray*}}
\newcommand{\eear}{\end{eqnarray*}}
\newcommand{\bt}{\begin{theorem}}
\newcommand{\et}{\end{theorem}}

\newcommand{\D}{\mathbb{D}}
\newcommand{\Dcal}{\mathcal{D}}
\newcommand{\Acal}{\mathcal{A}}
\newcommand{\Bcal}{\mathcal{B}}
\newcommand{\Hcal}{\mathcal{H}}
\newcommand{\Tcal}{\mathcal{T}_g(S)}
\newcommand{\R}{\mathbb{R}}
\newcommand{\C}{\mathbb{C}}
\newcommand{\CH}{\mathbb{CH}^2}

\newcommand{\V}{\mathcal{V}}
\newcommand{\Pat}{\mathcal{P}}

\newcommand{\Cal}{\mathcal{C}}
\newcommand{\E}{\mathcal{E}}

\newcommand{\dz}{\frac{\partial}{\partial z}}
\newcommand{\dbz}{\frac{\partial}{\partial \overline z}}

\newcommand{\dV}{\frac{i}{2}g_X(z)dz\wedge d\overline{z}}
\newcommand{\dzb}{d\overline{z}}
\newcommand{\db}{\overline{\partial}}

\numberwithin{equation}{section}

\allowdisplaybreaks

\def\Xint#1{\mathchoice
  {\XXint\displaystyle\textstyle{#1}}%
  {\XXint\textstyle\scriptstyle{#1}}%
  {\XXint\scriptstyle\scriptscriptstyle{#1}}%
  {\XXint\scriptscriptstyle\scriptscriptstyle{#1}}%
  \!\int}
\def\XXint#1#2#3{{\setbox0=\hbox{$#1{#2#3}{\int}$}
    \vcenter{\hbox{$#2#3$}}\kern-.5\wd0}}

\def\dashint{\Xint-}

\makeatletter
\def\@citestyle{\m@th\upshape\mdseries}
\def\citeform#1{{\bfseries#1}}
\def\@cite#1#2{{%
  \@citestyle[\citeform{#1}\if@tempswa, #2\fi]}}
\@ifundefined{cite }{%
  \expandafter\let\csname cite \endcsname\cite
  \edef\cite{\@nx\protect\@xp\@nx\csname cite \endcsname}%
}{}
\makeatother
\begin{document}
\parskip1ex
\title[New Parametrization of the space of {\cmli}s in $\CH$]
{Donaldson Functional in Teichm\"uller Theory}

\author{Zheng Huang}
\address[Z. ~H.]{Department of Mathematics, The City University of New York, Staten Island, NY 10314, USA}
\address{The Graduate Center, The City University of New York, 365 Fifth Ave., New York, NY 10016, USA}
\email{zheng.huang@csi.cuny.edu}

\author{Marcello Lucia}
\address[M.~L.]{Department of Mathematics, The City University of New York, Staten Island, NY 10314, USA.}
\address{The Graduate Center, The City University of New York, 365 Fifth Ave., New York, NY 10016, USA}
\email{marcello.lucia@csi.cuny.edu}

\author{Gabriella Tarantello}
\address[G.~T.]{Dipartimento di Matematica, Universit\`a di Roma ``Tor Vergata", Via della Ricerca Scientifica, I-00133 Roma, ITALY.}
\email{tarantel@mat.uniroma2.it}

\date{Nov. 20, 2021}


\begin{abstract}
In this paper we define a Donaldson type functional whose Euler-Lagrange equations are a system of differential equations which 
corresponds to Hitchin's self-duality equations for a suitable choice of Higgs bundle on closed {\RS}s. The main challenge of this functional is its 
lack of regularity and lack of compactness when defined in its natural domain of definition. Though a standard variational approach cannot directly 
be applied, we provide the appropriate analytical tools that make Donaldson functional treatable by a variational viewpoint. We prove that this 
functional admits a unique {\cp} corresponding to its global minimum. As an immediate consequence, we find that this system of self-duality equations 
admits a unique solution. Among the applications in geometry of this fact, we obtain a parametrization of closed constant mean curvature immersions in 
hyperbolic manifolds (possibly incomplete), and their moduli spaces. 
\end{abstract}

\maketitle

%\tableofcontents
%%%%%%%%%%%%%%%%%%%%%%%%%%%%%%%
\section {Introduction}
%\subsection{Geometric Background}
 
In his paper we present a variational approach to construct {\mi}s and {\mli}s of closed surfaces in {\htm}s and complex hyperbolic $2$-manifolds. 
In this way we obtain useful information about the representations of the fundamental group of the surface into $\text{PSL}(2,\C)$ and $\text{PU}(2,1)$.

Throughout this paper, we let $S$ be a smooth, closed, oriented surface of genus $g \ge 2$, and $\pi_1(S)$ be its {\fg}. The {\TS} 
$\mathcal{T}_g(S)$ is the space of {\cs}s on $S$, modulo biholomorphisms in the homotopy class of the identity. Uhlenbeck (\cite{Uhl83}) initiated a 
study of {\MS}s of {\mi}s of a closed surface into a {\tm} of constant sectional curvature~$-1$. Typically in this context, the {\tm} is hyperbolic, 
homeomorphic to $S\times\R$ and possibly incomplete. These {\mi}s naturally induce representations of $\pi_1(S)$ into the group $\text{PSL}(2,\C)$, the 
(orientation preserving) isometry group of $\mathbb{H}^3$. She considered the possibility of characterizing such class of irreducible representations, 
by fixing a conformal class $X$ on $S$ and a {\hqd} $q(z)dz^2$ on $X$. She proved a range of results for {\mi}s in {\qf} manifolds by using a bifurcation analysis 
based on the implicit function theorem. More recently, it was shown (\cite{HL12, HLT21}) that for a given data $(X,q(z)dz^2)$, a {\mi} may not exist. 
When it exists, one obtains a second solution in addition to the stable one constructed by Uhlenbeck. A similar multiplicity result was pointed out in 
\cite{HLL13} for prescribed {\hcd}s in the construction of {\mli}s. It is natural then to ask if one can parametrize the space of {\mi}s in another way.

In this spirit, Gon{\c{c}}alves-Uhlenbeck (\cite{GU07}) proposed to parametrize the space of immersions of {\cmc} surfaces in {\htm}s (homeomorphic 
to $S\times\R$ and possibly incomplete) by elements of the tangent bundle of the {\TS}. Such tangent bundle is identified as the collection of pairs 
$(X,[\beta])$, where $X$ is a {\cs} on the surface and $[\beta]$ is a {\coc} of $(0,1)$-forms valued in $T^{1,0}_X$, the holomorphic tangent bundle 
over $X$. We can trace back an analogous point of view to the Higgs bundle approach introduced by Hitchin in \cite{Hit87}. For example, we see that, 
for given a pair $(X,[\beta])$ one can obtain a {\mi} of $X$ in a {\htm} by solving the Gauss-Codazzi equations:
\be\label{GsCdz}	
   \left\{
   \begin{aligned}
     \partial\db\log(h)-h^2\|\beta\|^2-h^2=0\ ,\\
      \db(\ast_h\beta)=0\ 
   \end{aligned}
   \right.
\ene
expressed in terms of an Hermitian metric $h$ defined on the line bundle $K^{-\frac12}$, where $K= (T_X^{1,0})^{-1}$ is the canonical bundle of $X$, and a suitable 
representative in the class $[\beta]$ which, abusing our notations, is still denoted by $\beta$ in \eqref{GsCdz}. As usual we use the Hermitian extension of 
the {\hym} on $X$ to define $\|\beta\|^2$, while $\ast_h$ denotes the Hodge dual operator {\wrt} the metric $h$. From a solution of  \eqref{GsCdz} we 
obtain the pullback metric $g$ of the immersion from the Hermitian metric $h^2$ on $T^{1,0}_X$, and also we find that $\beta$ is harmonic {\wrt} $g$. 
Furthermore, $4\ast_h\beta$ defines a {\hqd} on $X$ whose real part identifies the {\sff} of the immersion.

We can show the equivalence between \eqref{GsCdz} and Hitchin's self-duality equations for a suitable nilpotent $\text{SL}(2,\C)$-Higgs bundle 
$(\E,\phi)$ given as follows. We let the rank two bundle $\E = K^{-\frac12}\oplus K^{\frac12}$ equipped with the holomorphic structure 
$\db_\E = \db + \begin{pmatrix}
0 & 0 \\
\beta & 0
\end{pmatrix}$ and Higgs field $\phi = \begin{pmatrix}
0 & 0 \\
1 & 0
\end{pmatrix}$. Then as explicitly derived in \cite{ALS21}, the pair $(h, \beta)$ satisfies \eqref{GsCdz} if and only if the Hermitian metric 
$H = \begin{pmatrix}
h& 0 \\
0 & \frac{1}{h}
\end{pmatrix}$ on $\E$ is the unique solution to Hitchin's self-duality equation:
\be\label{hit0}
F_{\nabla_H} + [\phi,\phi^{\ast_H}] = 0,
\ene
where $F_{\nabla_H}$ is the curvature form of the Chern connection ${\nabla_H}$, and $\phi^{\ast_H}$ is the Hermitian adjoint of 
$\phi$ {\wrt} $H$. In particular, by the general results of \cite{Hit87, Don87}, the given Higgs bundle $(\E,\phi)$ is stable. For full details of the one-to-one 
correspondence between the system \eqref{GsCdz} and Hitchin's self-duality equation \eqref{hit0}, we refer to \cite{ALS21}. See also recent survey articles 
(\cite{Wen16, Li19}) on Higgs bundles in relation to harmonic maps and other topics.

It turns out that problem \eqref{GsCdz} has a natural variational structure and the associated functional is referred by \cite{GU07} as {\df}, inspired by those introduced in 
K\"ahler geometry in (\cite{Don85, Sim88}) for the construction of Hermitian-Einstein metrics in holomorphic vector bundles.

We aim to pursue such a variational approach proposed in \cite{GU07}, and analyze the corresponding functional $\Dcal$ (see definition 1.5) for more general 
{\coc}es $[\beta]$ (dual to holomorphic $k$-differentials, for any $k\ge 2$), and prove a general uniqueness result showing that the global minimum is the only 
{\cp} of $\Dcal$. 
%Again this shows an one-to-one correspondence between the moduli space of the solutions for the Euler-Lagrange equations and the pair $(X,[\beta])$. 

Actually from (\cite{LM13}), we know that when $k=3$ there is a close relation between {\hcd}s on $X$ and equivariant {\mli}s from $\D$ into $\CH$. Our results provide 
an alternative proof for the one-to-one correspondence established in \cite{LM19} by means of Higgs bundles theory. On the other hand, for $k\ge 4$, a Higgs bundles theory 
approach to the existence and uniqueness issue is not available at present. It will be interesting and useful to interpret our solution in terms of a Hermitian metric solving 
Hitchin's equation for an appropriate choice of a ``stable" Higgs bundle.

More importantly, the variational approach provides us with the analytical framework to investigate the asymptotic behavior of minima corresponding 
to the data $(X,[t\beta])$ as $t \to \infty$, and this will be the main objective of our future work.

%\subsection{Donaldson Functional}

In the present paper we consider a Donaldson type functional defined in terms of the data $(X, [\beta])$, where $X$ is a {\cs} on the surface and $[\beta]$ is 
a {\coc} of $(0,1)$-forms valued in the bundle $E = \bigotimes\limits^{k-1} T^{1,0}_X \ (k\ge 2)$, namely $[\beta]\in \Hcal^{0,1}(X, E)$. To this purpose, 
let $\db$ be the induced holomorphic structure on $E$, and $g_X$ be the unique {\hym} on $X$, with Hermitian extension $h_X$. We let $\beta_0$ in 
$[\beta]$ be the unique harmonic element {\wrt} the {\hym}, so that $[\beta] = [\beta_0+\db\eta]$ for some section $\eta$ of $E$.

As usual, to remain in the given conformal class, we let $h=e^uh_X$ so that $g = e^{2u}g_X$. In this way, system \eqref{GsCdz} can be formulated in terms 
of the unknowns $(u,\eta)$ as follows:
\be\label{z}	
   \left\{
   \begin{gathered}
     \frac{ \Delta u+1}{4}-e^{2u}-\|\beta_0+\db\eta\|^2e^{2(k-1)u}=0\ \ \ \ \ {\text{ on }} X,\\
      \db\left(e^{2(k-1)u}\ast_E(\beta_0+\db\eta)\right)=0\ ,
   \end{gathered}
   \right.
\ene
where the Laplacian $\Delta$, the Hodge dual $\ast_E$, and the norm $\|\cdot\|$ are taken {\wrt} the background {\hym} $g_X$ and corresponding 
Hermitian extension $h_X$. To simplify notations it is convenient to operate the following change of variables: 
$u\mapsto 2(u+\ln2), \ \ [\beta]\mapsto [\frac{\beta}{2^{k-1}\sqrt{k-1}}]$, so that system \eqref{z} takes the form:

\be\label{hit}	
   \left\{
   \begin{gathered}
      \Delta u+2-2e^u-8(k-1)\|\beta_0+\db\eta\|^2e^{(k-1)u}=0\ \ \ \ \ {\text{ on }} X,\\
      \db\left(e^{(k-1)u}\ast_E(\beta_0+\db\eta)\right)=0\ .
   \end{gathered}
   \right.
\ene
Interestingly, solutions to system \eqref{hit} correspond to {\cp}s of the following functional:
\be\label{df}
\Dcal(u,\eta)=\int_X\{\frac14|\nabla u|^2- u+e^u+4\|\beta_0+\db\eta\|^2e^{(k-1)u}\}dA,
\ene
considered in its natural domain 
\be\label{domain}
\mathcal{W} = \{(u,\eta) \in H^{1}(X) \times W^{1,2}(X,E): \int_Xe^{(k-1)u}\|\beta_0+\db\eta\|^2dA < \infty\},
\ene
where $H^{1}(X)$ and $W^{1,2}(X,E)$ are the usual Sobolev spaces (see section 3 for details). In view of the connection with Hitchin's self-duality 
equations, as in \cite{GU07}, the functional $\Dcal$ is refered to as a ``Donaldson functional".

For $k=2$, the authors in \cite{GU07} observed through a formal computation that the second variation of $\Dcal$ at each possible {\cp} is positive 
definite. By this observation they claimed uniqueness of the {\cp}, as typically it would follow by standard global bifurcation arguments. However, 
$\Dcal$ is not continuous, or even weakly lower semicontinuous in $\mathcal{W}$, and therefore it needs particular care in order to be tackled by 
nonlinear techniques, as far as ``regularity" and ``compactness" issues are concerned.

Thus, to gain some ``regularity" for $\Dcal$, we work in the stronger Banach space $\V = H^{1}(X) \times W^{1,p}(X,E)$ with a fixed $p>2$. However, while 
$\Dcal$ is differentiable of any order in $\V$ (see Theorem ~\ref{smooth}), now we face a serious problem when verifying any sort of ``compactness" 
property for $\Dcal$ in $\V$, as for example the well known {\PSc}. For these reasons, the available variational approaches developed for nonsmooth 
functionals (e.g. the ``approximation" approach proposed by Struwe (\cite{Str08}) fail to apply to $\Dcal$. In fact, without ``compactness", even the 
information (we have obtained) that all {\cp}s of $\Dcal$ in $\V$ are strict local minima (see Proposition ~\ref{sphere}) isn't strong enough to imply 
``uniqueness". Indeed we could run into a situation similar to the function: $f(z) = |e^z-1|^2, \ z\in \C$, whose {\cp}s are infinitely many strict local 
minima exactly located at $z = 2\pi ni, \ n\in \mathbb{Z}$. Clearly the main difficulty is to gain control on the component $\eta$. However, such 
component is well behaved along {\cp}s, by the holomorphic condition in \eqref{hit}. We exploit exactly this fact, and by using Ekeland 
$\epsilon$-variational principle, we succeed to construct pre-compact {\PS} sequences. In this way we carry out a variational approach and show 
that indeed $\Dcal$ admits a global minimum in $\mathcal{W}$, which is its unique {\cp}.                   

\begin{thx}\label{main}
Let $X \in \Tcal$ be a closed {\RS}, and $E = \bigotimes\limits^{k-1} T^{1,0}_X$ be the tensor product of its holomorphic tangent bundles. Then for each 
{\coc} $[\beta]$ of $(0,1)$-forms valued in $E$, {\df} \eqref{df} admits a unique critical point $(u,\eta) \in \V$ corresponding to a global minimum. 
Furthermore, $(u,\eta)$ is smooth and it is the only solution to the system \eqref{hit}.
\end{thx}
\vskip 0.1in
There are several applications of Theorem A. For instance, in the case $k=2$ and $E=T^{1,0}_X$, the {\mi} provided by Theorem A can be lifted to a 
{\mi} from the universal cover $\D$ into $\mathbb{H}^3$ which is equivariant {\wrt} the associated representation $\rho: \pi_1(S)\rightarrow \text{PSL}(2,\C)$ 
(see Section 5 of \cite{Uhl83}). On the other hand, it is always possible to recover the pair of data $(X,[\beta])$ out of an equivariant {\mi} from $\D$ into 
$\mathbb{H}^3$ {\wrt} some representation. Recalling that a representations is irreducible if and only if the {\mi} is not totally geodesic (\cite{LM19b}), 
by virtue of Theorem A, we conclude: 
%establish  a one-to-one correspondence between the tangent bundle over $\Tcal$ (given by $\Tcal\times \Hcal^{0,1}(X, E)$) and the moduli space of equivariant {\mi}s, according to 
%the representation $\rho$ of $\pi_1(S)$ constructed by Uhlenbeck (see Corollary 5.2 of \cite{Uhl83}). In particular, we have  
\begin{corx}\label{maink2}
The moduli space of {\mi}s of $\D$ into $\mathbb{H}^3$ which are equivariant {\wrt} an irreducible representation of the fundamental group $\pi_1(S)$ 
into the group $\text{PSL}(2,\C)$ can be identified with the space $\Tcal\times \Hcal^{0,1}(X, E)\backslash\{0\}$.
\end{corx}
Corollary 1 was recently proved in (\cite{LM19b}) via a Higgs bundle approach (and they also attributed this to \cite{DEL97} from the point of view of 
birational algebraic geometry).

In case $k=3$, we have $E = T^{1,0}_X\bigotimes T^{1,0}_X$ and by Serre duality, $\Hcal^{0,1}(X, E)$ is isomorphic to the space $\Cal_3(X)$ of {\hcd}s on 
$X$. As before, these are used to parametrize the space of equivariant {\mli}s from $\D$ into complex hyperbolic plane $\CH$ (see for instance 
\cite{LM13, HLL13}). Our main theorem implies the following characterization of all such equivariant {\mli}s, as seen in \cite{LM19}: 
\begin{corx}\label{maink3}
The {\mli}s of $\D$ into $\CH$ which are equivariant {\wrt} an irreducible representation of $\pi_1(S)$ into the group $\text{PU}(2,1)$ 
are in one-to-one correspondence with the pairs $(X, [\beta]) \in \Tcal\times \Hcal^{0,1}(X, E)\backslash\{0\}$.
\end{corx}
Another application of Theorem A concerns the parametrization of the {\MS} of {\cmc} immersions of $S$ into germs of {\htm}s (see \cite{Tau04}). This 
problem can be reduced to study {\df} (with $k=2$) after a change of variable (see details in Section 6), and in view of Theorem A, it holds: 

%As mentioned in \S1.1, in \cite{GU07} the second variation of the functional was calculated at possible {\cp}s, and shown to be positive. But, as already mentioned, this information may not suffice to conclude uniqueness, unless we take good care of the ``compactness" issue. We prove:
\begin{corx}\label{cmc} 
For each given constant $c$ (with $c^2 < 1$), there is a one-to-one correspondence between the space of {\cmc} $c$ immersions in a germ of {\htm}s 
and the space $\Tcal\times \Hcal^{0,1}(X, E)$.
\end{corx}

The paper is organized as follows. In \S 2, after introducing the necessary notation, we focus on proving a Poincar\'e inequality with a sharp constant which is 
crucial to show that the Hessian of the functional $\Dcal$ at a {\cp} is positive definite.

We will break down the main result to the existence and uniqueness parts. In \S 3, we prove {\df} is well defined, bounded from below and smooth in the 
Banach space $\V$ defined above. We manage to establish appropriate ``compactness" properties for $\Dcal$ in $\V$ and construct a \underline{convergent} 
{\msq} yielding to the desired {\cp}. Furthermore, this minimum is regular by standard elliptic estimates.

Our strategy for proving the uniqueness part is done in two steps, contained in \S 4 and \S 5. We prove first that the quadratic form associated to the second 
variation of $\Dcal$ at a {\cp} $(u,\eta)$ is coercive in the space $H^1(X)\times W^{1,2}(X,E)$. Combining this fact with the holomorphic property of $\eta$ in 
\eqref{hit} and the convexity of the functional {\wrt} such variable, we succeed in establishing that actually the {\cp} $(u,\eta)$ is a strict local minimum for 
$\Dcal$ even {\wrt} the {\underline{stronger}} norm of the space $\V$ (Proposition ~\ref{sphere}). To conclude uniqueness we argue by contradiction, and by 
assuming that $\Dcal$ admits two distinct {\cp}s, by a ``{\mpa}" construction and a suitable use of Ekeland's $\epsilon$-variational principle, we produce an 
extra {\cp} for $\Dcal$ which can not be a strict local minimum (Theorem ~\ref{add}). This provides a full proof of our main result Theorem A.

%------------------------------------------------------
\subsection*{Acknowledgements}
We are deeply grateful for many helpful discussions with Qiongling Li. We also thank Sagun Chanillo and Song Dai and the referees for helpful suggestions. 
M.L. was supported by MINECO grant MTM2017-84214-C2-1-P. The research of G.T. was supported by MIUR excellence project, Department of Mathematics, 
University of Rome ``Tor Vergata" CUP E83C18000100006, and by ``Beyond Borders" research Project ``Variational approaches to PDE". 

%%%%%%%%%%%%%%%%%%%%%%%%%%%%%%%%%%%%%%%
\section{Background Material and Poincar\'e inequalities}
%We will collect some well-known results in Hodge theory, Hermitian geometry and we derive some technical lemmas which we will use later in the proof.
\subsection{Holomorphic Differentials and First Dolbeault Group}

% Since the paper will involve several different fields, 
In this subsection, we present the duality between the space of holomorphic $k$-differentials and the $(0,1)$-Dolbeault cohomology group. For this purpose, 
we collect some notations.
%Since there are no nonconstant holomorphic functions defined on $S$, holomorphic $k$-differentials ($k\ge 2$) 
%on {\RS}s and their duals have been essential parts of {\RS} theory, {\Tt} and the representations of surface groups in various character varieties. In 
%particular, when $k=2$, the space of {\hqd}s is identified with the cotangent space of {\TS} at a given {\cs}, while the dual space, the space of {\coc} of 
%so-called Beltrami differentials ($(-1,1)$-forms) is the tangent space (\cite{Ber72}). There are also deep connections with {\hqd}s and harmonic maps and 
%minimal surface theory (see for instance \cite{Uhl83}). When $k=3$, {\hcd}s appear in connection with {\mli} of a {\RS} (as a complex manifold) into the 
%complex hyperbolic plane $\mathbb{CH}^2$ (see for instance \cite{LM13, HLL13, LM19}). They also play a crucial role in understanding the 
%$SL(3,\R)$-{\cv} as they naturally appear in connection with real convex projective structures on a closed surface (\cite{Lof01, Lof07, Lab07}). In more 
%generality, Hitchin (\cite{Hit92}) studied the $SL(k,\R)$-{\cv} using holomorphic $k$-differentials. 

\ben[{\bf (1)}]
\item 
If $(x,y)$ is a real local coordinate where $z=x+iy$, then $\dz=\frac12(\frac{\partial}{\partial x}-i\frac{\partial}{\partial y})$ and 
$\dbz=\frac12(\frac{\partial}{\partial x}+i\frac{\partial}{\partial y})$, also $dz=dx+idy$ and  $d\overline{z}=dx-idy$;
\item
$X \in \Tcal$ is a {\cs} on $S$. It has a unique {\hym} denoted by $g_X$, whose volume form in a holomorphic coordinate $\{z\}$ can be written as 
$dA =\dV$, where $dz\wedge d\overline{z}=\frac12\{dz\otimes \dzb-\dzb\otimes dz\}$;
\item
Throughout this paper, for Laplace operator, inner products, norms and volume elements, we always use the {\hym} as the 
background metric unless we specify a lower index such as $\Delta_g$, $\langle\cdot,\cdot\rangle_g$, $\|\cdot\|_g$ or $dA_g$.
\item
$T^{1,0}_X$: the holomorphic tangent bundle over $X$. Since the complex dimension of $X$ is 1, the dual of $T^{1,0}_X$ coincides with the 
canonical bundle $K_X = K$.
\item
We always assume $k\ge 2$, and denote the tensor product $E =E_k = \bigotimes\limits^{k-1} T^{1,0}_X$, and $E^*$ the dual bundle of $E$;
%, and $\E$ the sheaf of holomorphic sections of the bundle $E$;
\item
$A^{0}(E) = \{\eta\}$ is the space of smooth sections of $E$; we denote $A^{0,1}(X,\C)$ the space of complex-valued $(0,1)$-forms on $X$, and 
$A^{0,1}(X,E) = \{\beta\}$ the space of $(0,1)$-forms on $X$ valued in $E$, i.e., $A^{0,1}(X,E) = A^{0,1}(X,\C)\otimes E$;
\item
$\Cal_k(X)$: the space of {\hkd}s on $X$, or equivalently the space of holomorphic sections of the bundle $\bigotimes\limits^{k}(T^{1,0}_X)^*$, often also  
denoted by $H^0(X,\bigotimes\limits^{k} (T^{1,0}_X)^*)$. Any such differential is locally of the form $q(z)dz^k$ on $X$, where $q(z)$ is holomorphic. 
As a consequence of the Riemann-Roch Theorem, the complex dimension of $\Cal_k(X)$ is $(2k-1)(g-1)$; %or equivalently, $(1,0)$-form valued in $E^*$, 
\item
The $(0,1)$-Dolbeault cohomology group $\Hcal^{0,1}(X, E)$ is defined as the quotient space $\faktor{A^{0,1} (X,E)}{\db(A^{0}(E))}$, where 
$\db: A^{0}(E) \to A^{0,1}(X,E)$ is the d-bar operator. By Hodge Theory, there is a natural isomorphism between $\Hcal^{0,1}(X, E)$ and the space of 
harmonic $(0,1)$-forms valued in $E$. Therefore $\forall [\beta] \in \Hcal^{0,1}(X, E)$, there is a unique harmonic element $\beta_0$ {\wrt} the {\hym} 
such that $\beta = \beta_0+\db\eta$, with $\eta \in A^0(E)$.
%There is a natural Serre duality between $\Cal_k(X)$ and $\Hcal^{0,1}(X, E)$ which we describe in \S 2.1.
\een
%We use the following notations:\ben[\bf (i)]\item$h(w,v)=g(w,\overline{v})$ is the corresponding {\hem}, and $|w|^{2}=h(w,w)$.\een

%We define the {\Rm} $g$ on forms similarly as the wedge products of orthonormal basis. Then we can define {\hem} $h$ on forms in the same 
%manner as above. Note that the metric on forms does not coincide with the metric on tensors if we write forms out as alternate sums of tensors. 
%By abusing the notations, we denote %\beq |dz|^{2}=|\dzb|^{2}=g^{-1}, \quad |dz\wedge \dzb|^{2}=g^{-2}.\eeq
Given a Riemannian metric $g$ on $X$, it induces a {\hem} $h$ on $X$ with $h(v,w)=g_{\C}(v,\overline{w})$ where $g_{\C}$ is the complex extension 
of $g$ as a bilinear form. We can extend it to obtain a {\hem} $\langle\cdot,\cdot\rangle_g$ on sections and forms valued in $E$.

Let $\alpha$ be a differential form on $X$ valued in $\C$, then we define Hodge star $*\alpha$ as $\psi\wedge*\alpha=\langle\psi,\alpha\rangle dA$ for 
every $\psi$. In a local holomorphic coordinate $\{z\}$, we have $\ast dz=i\dzb$ and $\ast\dzb=-idz$. This can be extended to differential forms valued 
in $E$ as follows. We define $\flat$ operator to be the identification of $E$ and $E^{*}$, i.e. for $w,v\in A^0(E)$, $(\flat w)(v)=\langle v,w\rangle$. 
We also define the metric on $E$-valued forms as the composition of the two metrics on forms and $E$. Now we define the Hodge star operator 
$*_E: A^{0,1}(X,E) \to A^{1,0}(X,E^*)$ as the conjugate linear map $*_E(\alpha_0\otimes e) = *\alpha_0\otimes\flat e$, with 
$\alpha_0 \in A^{0,1}(X,\C)$ and $e\in E$. We can also define wedge product for forms valued in vector bundles. Particularly, if 
$\alpha = \alpha_0\otimes e^\ast \in A^{1,0}(X,E^\ast)$ and $\beta = \beta_0\otimes e \in A^{0,1}(X,E)$, with $\alpha_0 \in A^{1,0}(X,\C)$ and 
$\beta_0 \in A^{0,1}(X,\C)$, $e \in A^0(E)$ and $e^\ast \in A^0(E^\ast)$, then  
\beq
\alpha\wedge\beta = e^\ast(e)\alpha_0\wedge\beta_0. 
\eeq
With this and the definition of $\ast_E$ on forms valued in bundles, for $\beta_1, \beta_2 \in A^{0,1}(X,E)$, we note:
\be\label{well}
\ast_E\beta_1\wedge\beta_2= \langle\beta_1,\beta_2\rangle dA.
\ene

%Let us work with the duality in a slightly more general setting. 
We have the following natural bilinear form:
\beq
A^{1,0}(X,E^\ast)\times A^{0,1}(X,E) \rightarrow \C, \ (\alpha, \beta)\mapsto \int_X\alpha\wedge\beta.
\eeq
Since we have the identification $A^{1,0}(X,E^\ast)\cong A^{0}(X,E^*\otimes K_X)$, and noting that Stoke's Theorem gives 
$\int_X\alpha\wedge\db\eta = 0$ for all $\alpha \in A^{1,0}(X,E^\ast)$ which are holomorphic, and all $\eta\in A^{0}(E)$, we obtain the 
following bilinear form:
%\beqH^{0}(X,E^*\otimes K_X)\times A^{0,1}(X,E) \rightarrow \C, \ (\alpha, \beta)\mapsto \int_X\alpha\wedge\beta.\eeq
%where $K_X$ is the canonical line bundle, which in our case is the $(1,0)$-part of the complexified cotangent bundle. 
%This induces a bilinear form:
\beq
H^{0}(X,E^*\otimes K_X)\times \Hcal^{0,1}(X,E) \rightarrow \C, \ (\alpha, [\beta])\mapsto \int_X\alpha\wedge\beta.
\eeq
%Note that $X$ has complex dimension $1$. 
%This pairing \eqref{pairing} can be computed in the Dolbeault cohomology as follows: if $\alpha_1\in A^{0,1}(X,E), \alpha_2\in A^{0,0}(K_X\otimes E^*)$, 
%we define the pairing
%\be\langle\alpha_1,\alpha_2\rangle=\int_X \alpha_{1x}\wedge\alpha_{2x},\ene
%where the differential form $\alpha_{1x}\wedge\alpha_{2x}$ of degree $2$ is obtained by taking the exterior
%product on the forms, contracting $E$ and $K_X\otimes E^*$. One checks by Stokes' theorem that
%\beq\langle\alpha_1,\overline{\partial}\alpha_2\rangle=\langle\overline{\partial}\alpha_1,\alpha_2\rangle.\eeq
%Therefore if $\alpha_1$ and $\alpha_2$ are $\overline{\partial}$-closed, then $\langle\alpha_1,\alpha_2\rangle$ depends only on their {\coc}es: it is the desired pairing 
%(\ref{pairing}).   
This is a nondegenerate bilinear form and therefore it induces 
%between the cohomology groups $H^1(X,E)$ and $H^{0}(X,E^*\otimes K_X)$ is perfect, namely, the operator $*_E$ gives 
an injective homomorphism between $H^{0}(X,E^*\otimes K_X)$ and $\Hcal^{0,1}(X,E)^\ast$, which by Serre's duality Theorem (see \cite{Voi07}), 
is an isomorphism. As a consequence, given $[\beta_0+\db\eta] \in \Hcal^{0,1}(X,E)$, where $\beta_0$ is harmonic, by considering the following linear map 
\beq
\Hcal^{0,1}(X,E) \rightarrow \C, \ [\xi]\mapsto \int_X\langle\xi,\beta_0\rangle dA
\eeq 
(well defined since $\int_X \langle\db\eta,\beta_0\rangle dA = 0$), there exists a unique $\tilde\alpha\in H^{0}(X,E^*\otimes K_X)$ such that 
$\int_X \tilde\alpha\wedge\xi = \int_X\langle\xi,\beta_0\rangle dA$, i.e., $\tilde\alpha=\ast_E\beta_0$. Since in our notation, $\Cal_k(X) = H^{0}(X,E^*\otimes K_X)$, 
we obtain the following isomorphism:
\be\label{pairing}
\Hcal^{0,1}(X,E) \to \Cal_k(X), \ [\beta]\mapsto \ast_E\beta_0.
\ene

%%%%%%%%%%%%%%%%%%%%%%%%%%
\subsection{Poincar\'e type inequalities}
In this subsection, we will present two Poincar\'e type inequalities in Proposition ~\ref{Bochner} and Proposition ~\ref{PI}. We suspect inequalities of this type are standard, 
and we only include the proofs here for the sake of completeness.
%From now on, we denote $\ast = \ast_E$. 
%This identity expresses the anti-holomorphic Laplace operator acting on $\Omega(X,E)$ in terms of its conjugate 
%and some extra terms involving the curvature of $E$ and the torsion of the {\hem} $\omega$ (in the case that $\omega$ is not K\"{a}hler).
%We now apply {\BKN} to our setting ($r=n=1$ and $E = \bigotimes\limits^{k-1} T^{1,0}_X$) to obtain the following lemma, which can be 
%viewed as a version of the $L^2$-Poincar\'e inequality for differential forms in our setting.

\bpo\label{Bochner}
Let $X$ be a closed {\RS} with the Riemannian metric $g$, and $K(g)$ be its Gaussian curvature. Then 
\be\label{ineq1}
\int_X\langle\db_E \psi,\db_E\psi\rangle_g dA_g\geq -(k-1)\int_X K(g)\langle\psi,\psi\rangle_g dA_g,
\ene
holds for every section $\psi \in W^{1,2}(X,E)$. In particular if we use the {\hym} $g_X$, and $dA$ its volume form, we have
\beq
\int_X\langle\db_E \psi,\db_E\psi\rangle dA\geq (k-1)\int_X \langle\psi,\psi\rangle dA
\eeq
\epo
\bp
By density of smooth sections in $W^{1,2}(X,E)$, it is enough to prove this for smooth section $\psi$. In the conformal coordinate $\{z\}$ of $X$, we 
write the metric $g=e^{2\phi}|dz|^2$ and let $\omega =\frac{i}{2}e^{2\phi}dz\wedge d\bar{z}$ be its {\kf}. Then we have $\ast\omega=1$ and $\ast1=\omega$.

The exterior differentiation $d$ does not extend to vector-valued forms, we will work with Chern connection on holomorphic bundles here. 
For a holomorphic vector bundle $E$ and a local holomorphic frame $F$, the Chern connection $D_E=D'+D''$ can be characterized {\wrt} $F$ 
as follows:
\beq
D'=\partial_E+\theta(F),\ \ \ \ \ D''=\db_E,
\eeq
where $\theta(F)$ is the connection matrix {\wrt} $F$.

Writing $\Delta''=D''\delta''+\delta'' D''$ and $\Delta'=D'\delta'+\delta'D'$, where $\delta'$ and $\delta''$ are the formal adjoint operators of $D'$ and $D''$, respectively, 
we have the Chern curvature form $\Theta(E)=d\theta+\theta\wedge\theta$. We define the operators 
\beq
0 \xrightarrow{L} A^0(X,\C) \xrightarrow{L} A^{1,1}(X,\C), 
\eeq 
where the second arrow is defined as $\eta\mapsto \omega\wedge\eta = \eta\omega$, and its adjoint 
\beq
\Lambda: A^{1,1}(X,\C)\to A^0(X,\C), \ \Lambda = \ast^{-1}L\ast.
\eeq 
We extend them on forms valued in the bundle $E$ by tensor product with the identity map on $E$.

Using the Akizuki-Nakano identity (\cite{AN54}): $\Delta''=\Delta'+[i\Theta(E),\Lambda]$, we obtain:
\beq
\int\langle\Delta'' \psi,\psi\rangle_g =\int\langle\Delta' \psi,\psi\rangle_g+\int\langle[i\Theta(E),\Lambda]\psi,\psi\rangle_g.
%&=&\int\langle D' \psi,D' \psi\rangle_g+\int \langle[i\Theta(E),\Lambda]\psi,\psi\rangle_g.
\eeq
For a smooth section $\psi$ of $E$, we have $\delta'(\psi)=\delta''(\psi)=0$, and $D''=\db_E$, therefore,
\bear
\int\langle\db_E\psi,\db_E\psi\rangle_g&=&\int\langle D' \psi,D' \psi\rangle_g+\int \langle[i\Theta(E),\Lambda]\psi,\psi\rangle_g \\
&\ge&\int \langle[i\Theta(E),\Lambda]\psi,\psi\rangle_g.
%&=&\int\langle\Delta' \psi,\psi\rangle_g+\int 4(k-1)e^{-2\phi}\phi_{z\overline{z}}\langle\psi,\psi\rangle_g\\&=&\int\langle D' \psi,D' \psi\rangle_g
%+\int 4(k-1)e^{-2\phi}\phi_{z\overline{z}}\langle\psi,\psi\rangle_g,
%&=&\int\langle D' \psi,D' \psi\rangle_g-(k-1)\int K(g)\langle\psi,\psi\rangle_g, \text{~since~} K(g)=-4e^{-2\phi}\phi_{z\overline{z}}. \\
%&\geq& -(k-1)\int K(g)\langle\psi,\psi\rangle_g.
\eear
Let us calculate the term $[i\Theta(E),\Lambda]$ for our case where $E = \bigotimes\limits^{k-1} T^{1,0}_X$. Since the frame $F$ is locally 
$\bigotimes\limits^{k-1}\frac{\partial}{\partial z}$, we have $\theta(F)=\|F\|^{-2(k-1)}\partial \|F\|^{2(k-1)}=2(k-1)\partial \phi$, and 
$\Theta(E)=d\theta+\theta\wedge\theta=2(k-1)\db\partial \phi$, given that the complex dimension of the surface is one. 
%Therefore we have $\Delta''=\Delta'+[i2(k-1)\overline{\partial}\partial \phi,\Lambda]$.
Therefore, for each $\psi \in A^0(E)$, we get 
\beq
\frac{1}{2(k-1)}\Theta(E)(\psi) = (\db\partial \phi)(\psi) = \phi_{\bar{z}z}\psi d\bar{z}\wedge dz = -\phi_{z\bar{z}}\psi dz\wedge d\bar{z}.
\eeq
Thus, 
\beq
\frac{1}{2(k-1)}i(\Lambda\circ\Theta(E))(\psi) = -i\Lambda(\phi_{z\bar{z}}dz\wedge d\bar{z})\otimes\psi = -2\Lambda(\omega)\otimes (e^{-2\phi}\phi_{z\bar{z}})\psi.
\eeq
Since $\Lambda\psi = 0$, and $\Lambda(\omega)=\ast^{-1} L\ast \omega=\ast^{-1} L 1=\ast^{-1}\omega=1$, we have 
\beq
[i\Theta(E),\Lambda]\psi = i\Theta(E)\Lambda\psi - i\Lambda\Theta(E)\psi=4(k-1)e^{-2\phi}\phi_{z\bar{z}}\psi.
%&=&- i\Lambda(\Theta(E)(\psi) \\
%&=&2(\Lambda\frac12 \phi_{z\overline{z}}dz\wedge d\bar{z})\psi\\
%&=&2(\Lambda w\otimes e^{-2\phi}\phi_{z\bar{z}})\psi\\
\eeq
Using $K(g)=-4e^{-2\phi}\phi_{z\bar{z}}$, we find:
\bear
\int\langle\db_E\psi,\db_E\psi\rangle_g&\ge&\int 4(k-1)e^{-2\phi}\phi_{z\bar{z}}\langle\psi,\psi\rangle_g \\
&=& -(k-1)\int K(g)\langle\psi,\psi\rangle_g.
\eear
\ep
%Lemma ~\ref{Bochner} is simply saying we have:
%\beq \|\overline{\partial} \psi\|_{L^2}\geq \sqrt{k-1}\|\psi\|_{L^2}.\eeq

We end this section with the following $L^p$-version of the Poincar\'e inequality, and we state it for a general holomorphic vector bundle $\tilde{E}$ over $X$:
\bpo\label{PI}
Let $X$ be a closed {\RS} with the {\hym} $g_X$, and $\tilde{E}$ be a holomorphic vector bundle over $X$, and $h$ be an Hermitian metric on 
$\tilde{E}$. Suppose the only global holomorphic section of $\tilde{E}$ is zero. Then for $1<p<\infty$, 
there exists a positive constant $C=C(X,\tilde{E},h,p)$, such that, for any section $\psi$ of $W^{1,2}(X,\tilde{E})$, there holds: 
\be\label{Lp}
\int_X \|\db_E\psi\|^{p}dA\geq C\int_X \|\psi\|^{p}dA.
\ene
\epo
%Note that in our case, $E = \bigotimes\limits^{k-1} T^{1,0}_X$ is a negative bundle since the genus of the {\RS} $X$ is at least two, 
%which means zero section is the only global holomorphic section. 
\bp
We prove this by contradiction. If not, there exists a sequence of smooth sections $\psi_{j}$, such that
\begin{eqnarray}
||\psi_{j}||_{L^p}=1~\text{and}~||\db_E\psi_j||_{L^p}<\frac{1}{j},\quad\text{for}~j=1,2,\cdots.\label{s}
\end{eqnarray} 
We claim that, regarding $\db_E$ as a real operator, it is elliptic. Let $e_{\alpha}$ be a local holomorphic frame of $\tilde{E}$, then locally, we can write 
$\psi=\psi^{\alpha}e_{\alpha}$, where $\psi^{\alpha}$ is a complex valued function.

Let $\psi^{\alpha}=u^{\alpha}+iv^{\alpha}$ and ${z=x+iy}$ be a local holomorphic coordinate system on $X$. Then $d\bar{z}\otimes e_{\alpha}$ 
is a local frame for $(T^{1,0}_X)^\ast\otimes \tilde{E}$, and locally, we have 
\bear
\db_E\psi&=&\frac{1}{2}(\partial_{x}(u^{\alpha}+iv^{\alpha})-i\partial_y(u^{\alpha}+iv^{\alpha}))d\bar{z}\otimes e_{\alpha}\\
&=&\frac{1}{2}((\partial_xu^{\alpha}+\partial_yv^{\alpha})+i(-\partial_yu^{\alpha}+\partial_xv^{\alpha}))d\bar{z}\otimes e_{\alpha}.
\eear
Considering now $\db_E$ as a real differential operator, we have 
\begin{displaymath}
\db_E:
\left(
\begin{array}{c}
u^{\alpha}\\
v^{\alpha}
\end{array}
\right)
\mapsto
\frac12\left(
\begin{array}{cc}
\xi_x & \xi_y\\
-\xi_y & \xi_x
\end{array}
\right)
\left(
\begin{array}{c}
u^{\alpha}\\
v^{\alpha}
\end{array}
\right)d\bar{z}\otimes e_{\alpha}.
\end{displaymath}
Since the matrix symbol associated to each component of the operator is given by 
$
\frac12\left(
\begin{array}{cc}
\xi_x & \xi_y\\
-\xi_y & \xi_x
\end{array}
\right),$
and it has nonzero determinant, we deduce that $\db_E$ is elliptic. (Notice that this is not true in general for high dimensional base manifolds). 
Then by standard elliptic estimates (\cite{DN55}), there exists a positive constant $C$, independent of $j$, such that
\beq
\|\psi_j\|_{W^{1,p}}\leq C(\|\db_E\psi_j\|_{L^{p}}+\|\psi_j\|_{L^{p}}).
\eeq
Then by (\ref{s}), we see that the $W^{1,p}$-norm of $\psi_j$ is uniformly bounded. Then we can choose a subsequence which converges weakly in 
$W^{1,p}$. We denote the weak limit by $\psi_0$. Then from the compact embedding theorem, $\psi_j$ converges strongly in $L^p$ to $\psi_0$. This 
implies $\|\psi_0\|_{L^{p}}=1$. On the other hand, since the functional $\int_X \|\db_E\psi\|^p dA$ is convex, then by weak lower semi-continuity,
\beq
\int_X \|\db_E\psi_0\|^p dA\leq \liminf \int_X \|\db_E\psi_j\|^p dA=0.
\eeq
Hence $\db_E\psi_0=0$, almost everywhere. From the elliptic regularity (see for instance \cite{Nar92}), $\psi_0$ is smooth, then holomorphic. From 
the assumption that there is no global holomorphic section, $\psi_0$ must be $0$, which contradicts the assumption that $\|\psi_0\|_{L^{p}}=1$. 
This completes the proof. 
\ep
\br
Proposition ~\ref{Bochner} implies that the only global holomorphic section for the negative holomorphic line bundle $E$ is zero. This fact can 
also be seen from standard Kodaira-vanishing theorems in complex geometry (see [Chapter VI, Theorem 2.4(b), \cite{Wel08}]). Therefore 
the assumptions in Proposition ~\ref{PI} hold for the bundle $E$. In the case of $L^2$-version of the inequality, the constant is 
explicit and universal, while the constant for general $L^p$-version is less explicit. We later will take advantage of both 
properties. %Proposition ~\ref{Bochner} is used in the proof of Lemma ~\ref{Poi}, and Proposition ~\ref{PI} is used in the proof of Proposition ~\ref{regu}.
\er
%%%%%%%%%%%%%%%%%%%%%%%%%%%%
%%%%%%%%%%%%%%%%%%%%%%%
\section{Existence of critical point for the Donaldson Functional}
Recall that the Sobolev space $H^1(X)$ is the closure of $C^1$-functions on {\RS} $X$, {\wrt} the $H^1$-norm defined as:
\beq
\|u\|_{H^1} = \big\{\int_X (|u|^2+|\nabla u|^2)dA\big\}^{\frac12}.
\eeq
Let us also define some Sobolev spaces for differential forms valued in $E$. The Hermitian form $\langle\cdot,\cdot\rangle$ defined fiberwise on each space 
$A^{0,j}(X,E)$ ($j=0,1$) induces a norm denoted by $\|\cdot\|$. In local coordinates, if $\eta\in A^{0,0}(X,E)=A^0(E)$, and $\beta \in A^{0,1}(X,E)$, then 
\beq
\|\eta\| = |\eta(z)|(g_X)^{\frac{k-1}{2}}, \ \ \ \  \|\beta\| = |\beta(z)|(g_X)^{\frac{k-2}{2}}.
\eeq
Given $q \ge 1$, a $(0,j)$-form $\alpha$ valued in $E$ is said to belong to $L^q(A^{0,j}(X,E))$ if $\int_X\|\alpha\|^qdA < \infty$. We also say a section $\eta$ of the 
bundle $E$ belongs to $W^{1,q}(X,E)$ if
%\beq
$\int_X(\|\eta\|^q + \|\db\eta\|^q)\ dA< \infty$,
%\eeq
and its $W^{1,q}$-norm is given as follows:
\be\label{w1p}
\|\eta\|_{W^{1,q}(X,E)}=\big\{\int_X(\|\eta\|^q + \|\db\eta\|^q)dA\big\}^{\frac1q}.
\ene
We analyze the functional $\Dcal$ in the space $\V= H^{1}(X) \times W^{1,p}(X,E)$, for fixed $p>2$, endowed with the norm $\|\cdot\|_{\V}$ given as follows
\be\label{vnorm}
\|(u,\eta)\|_{\V} = \sqrt{\|u\|_{H^1}^2+\|\eta\|_{W^{1,p}(X,E)}^2}.
\ene

We observe that $(\V, \|\cdot\|_{\V})$ is a uniformly convex Banach space (\cite{Cla36}), a property we explore in section \S5. We will prove 
the existence and regularity of a {\cp} for $\Dcal$.

The fact that $\Dcal$ is well defined is a consequence of the Moser-Trudinger inequality, which states that on any closed Riemannian surface $(X,g_X)$, 
there exists some constant $C_X > 0$ such that (see Theorem 2.50 in \cite{Aub98} or \cite{Fon93}):
\be\label{MTI}
\int_X e^{4\pi(\frac{u}{\|\nabla u\|_{L^2}})^2}dA \le C_X, \ \ \forall u\in H^1(X) \text{ with } \int_X udA = 0.
\ene  
This inequality has the following consequence:
\be\label{MTI1}
\int_X e^vdA \le C_Xe^{\dashint_X vdA}e^{\frac{\|\nabla v\|^2_{L^2}}{16\pi}}, \ \ \forall v\in H^1(X).
\ene
Since $\big|e^v-e^{v_0} \big|\le e^{v_0}e^{|v-v_0|}|v-v_0|$ (by mean value theorem) for any $v, v_0 \in H^1(X)$, using \eqref{MTI1}, we have:
\be\label{MTI2}
\int_X\big|e^v-e^{v_0} \big|dA \le \big(C_X\int_Xe^{4v_0}dA\big)^\frac14e^{\{\frac{\|v-v_0\|_{L^1}}{|X|}+ \frac{1}{4\pi}\|\nabla (v-v_0)\|^2_{L^2}\}}\|v-v_0\|_{L^2}.
\ene
In particular, these inequalities imply that, for every $q\ge 1$, the map $H^1(X)\to L^q(X), \ v \mapsto e^v$ is well defined, $C^\infty$, and compact.
%Though natural space to work in for $\Dcal$ is $H^{1}(X) \times W^{1,a}(X,E)$ for some 
%$a \in [1,2)$ (~\ref{apriori}), our strategy is to work in the space $\V$. We prove that $\Dcal$ is differentiable in $\V$ (Theorem ~\ref{smooth}), and 
%then we use a map $\Dcal_u(\eta)$ (\S 3.3), for each {\msq} in $H^{1}(X) \times W^{1,a}(X,E)$, to produce a {\msq} in $\V$. We will also prove 
%that One of the key geometric ingredients we use is that the space of {\hkd}s on a closed {\RS} is of finite dimensional.

%%%%%%%%%%%%%%%%%%%%%%%%%%%%%%%
\subsection{Differentiability of {\df}}
In this subsection we show that $\Dcal$ is smooth in the space $\V$. This justifies why we choose to work with such space, even though it does introduce 
some difficulties when dealing with ``compactness" issues. To this purpose, we first notice that {\df} 
\beq
\Dcal(u,\eta)=\int_X\big\{\frac14|\nabla u|^2- u+e^u+4\|\beta_0+\db\eta\|^2e^{(k-1)u}\big\}dA
\eeq
can be written as  
$\Dcal (u, \eta) = \Acal(u) + 4\Bcal (u, \eta)$, where
\be\label{D1}
\Acal(u) = \int_X \big\{\frac14|\nabla u|^2- u+e^{u}\big\}dA,
\ene
and
\be\label{D2}
\Bcal(u, \eta) = \int_X \|\beta_0+\db\eta\|^2e^{(k-1)u}dA.
\ene
It is well known that the functional $\Acal$ is well defined and $C^\infty$ in $H^1(X)$, and in particular, at each $u \in H^1(X)$, we have,  
\be\label{Aprime}
\Acal'_{u}[v] = \int_X\big\{\frac12 \nabla u\nabla v -v +e^uv\big\}dA, \ \forall v\in H^1(X),
 \ene
 and
 \be\label{Apprime}
\Acal''_{u}[v_1, v_2] = \int_X\big\{\frac12 \nabla v_1\nabla v_2 +e^uv_1v_2\big\}dA, \ \forall v_1, v_2\in H^1(X).
 \ene
 To show the regularity of $\Bcal$, we need the following lemma:
 \bl\label{mapLp}
 \ben
 \item
 The map $\tau: L^p(A^{0,1}(X,E)) \longrightarrow L^{\frac{p}{2}}(X),\ \ F \longmapsto \langle F, F\rangle$ is $C^\infty$ and we have: 
 \beq
 \tau'_F[\zeta] = 2\Re \langle F, \zeta\rangle, \ \  \tau''_F[\zeta_1,\zeta_2] = 2\Re \langle\zeta_1, \zeta_2\rangle,
 \eeq
 with vanishing higher derivatives;
\item
The map $\tilde\tau: L^p(A^{0,1}(X,E))\times L^{\frac{p}{p-2}}(X) \to L^1(X), \  (F,g) \mapsto \langle F, F\rangle g$ is also 
$C^\infty$, and we have 
\beq
\tilde\tau'_{(F,g)}[\zeta,v] = v\langle F, F\rangle + 2g\Re\langle F, \zeta\rangle,
\eeq
and 
\beq
\tilde\tau''_{(F,g)}[(\zeta_1,v_1), (\zeta_2, v_2)] = 2\Big(v_1\Re\langle F, \zeta_2\rangle+v_2\Re\langle F, \zeta_1\rangle
+ g\Re\langle \zeta_1, \zeta_2\rangle\Big).
\eeq
\een
 \el
 \bp
(i) The expressions of $ \tau'_f[\zeta]$ and $\tau''_f[\zeta_1,\zeta_2] $ are obtained by explicit calculations. Furthermore, setting 
$Y = L^p(A^{0,1}(X,E))$ and $Z = L^{\frac{p}{2}}(X)$, and using Cauchy-Schwarz inequality, we check the continuity of the following maps 
$$
   \begin{array}{rl}
       Y &\to  \mathcal{L} (Y, Z),       \\
       F &\mapsto   \tau'_F 
   \end{array}
   \quad
    \begin{array}{rl}
       Y &\to  Bil(Y \times Y ,Z),       \\
       F &\mapsto   \tau''_F 
   \end{array}
$$
where $ \mathcal{L} (Y, Z)$ is the space of continuous linear maps from $Y$ to $Z$, and $Bil(Y \times Y ,Z)$ is the space of continuous 
bilinear maps from $Y\times Y$ to $Z$, both endowed with the sup-norm.
%\beq |\langle f_1, f_2\rangle_E|\le\|f_1\|\|f_2\|.\eeq
\vskip 0.1in
\noindent
(ii) The map $\tilde\tau$ is a composition of the map 
  \begin{align*}
  \tau\times id: L^p(A^{0,1}(X,E))\times L^{\frac{p}{p-2}}(X) \longrightarrow L^{\frac{p}{2}}(X)\times L^{\frac{p}{p-2}}(X) &\\
  (F,f) \longmapsto (\langle F, F\rangle, f)
   \end{align*}
with the map $L^{\frac{p}{2}}(X)\times L^{\frac{p}{p-2}}(X)  \longrightarrow L^1(X), \  (g_1,g_2) \longmapsto g_1g_2$.
 \ep
\bt\label{smooth}
The functional $\Bcal \in C^\infty(\V)$, and consequently $\Dcal$ is smooth in $\V$. Moreover, let $\beta = \beta_0+\db\eta$,
 then for every $(u,\eta) \in \V$, we have  
\be\label{frechet}
\Bcal'_{(u,\eta)}[v,\ell] = \int_X e^{(k-1)u}\big((k-1)\|\beta\|^2v+2\Re\langle\beta,\db\ell\rangle \big), \ \forall \ (v,\ell) \in \V,
\ene
and 
\begin{eqnarray}\label{frechet2}
\Bcal''_{(u,\eta)}\langle[v_1,\ell_1], [v_2,\ell_2]\rangle &=& \int_X e^{(k-1)u}\big\{(k-1)^2\|\beta\|^2v_1v_2+ 2\langle\db\ell_1,\db\ell_2\rangle \notag\\ 
&& \ \ \ \ \ \ \ \  + 2(k-1)(v_1\Re\langle\beta,\db\ell_2\rangle + v_2\Re\langle\beta,\db\ell_1\rangle)\big\}, 
\end{eqnarray}
for all $(v_1,\ell_1),\ (v_2,\ell_2) \in \V$.
\et
\bp
From \eqref{MTI2}, we see that for each $q\ge 1$, the map $u \mapsto e^{(k-1)u}$ from $H^1(X)$ to $L^q(X)$ is smooth, and its derivative can 
be easily calculated. At this point we write $\Bcal(u,\eta) = \int_X \tilde\tau(\beta_0+\db\eta, e^{(k-1)u})dA$. Now formulas \eqref{frechet} and \eqref{frechet2} 
follow easily via direct calculations. 
\ep
%%%%%%%%%%%%%%%%%%%%%%%%%%%%%%%%%%%%%%%%%
\subsection{Characterization of {\cp}s} 
From now on, we use $* $ to denote the Hodge star $*_E$ for forms valued in $E$. The main result in this subsection is to characterize 
the {\cp}s of {\df}: 
\bt\label{EL}
The pair $(u,\eta) \in \V$ is a {\cp} for $\Dcal(u,\eta)$ if and only if $(u,\eta)$ is a \underline{smooth} solution for the system \eqref{hit}.
\et
\bp
For any family of $(u_t,\eta_t)$ with $\beta_t = \beta_0+ \eta_t$ and $u_t=u+tv,\eta_t=\eta+t\ell$, where $(v,\ell) \in \V$, we can readily 
compute the first variation of $\Dcal$ from \eqref{Aprime} and \eqref{frechet} as follows:
\bear
d/dt|_{t=0} \Dcal(u_t,\beta_t)&=&\int_X \Big\{\frac12 \nabla u\nabla v+(-1+e^u+4(k-1)\|\beta_0+\db\eta\|^2 e^{(k-1)u})v\\
&&\ \ \ \ \ \ \ \ \ \ \ \ \ +8\Re\langle\beta_0+\db\eta,\overline{\partial}\ell\rangle e^{(k-1)u}\Big\}\ dA
\eear
Hence we find that $(u,\eta)$ is a {\cp} of $\Dcal$ if and only if $u$ is a weak solution in $H^1(X)$ of the equation  
\be\label{cp-1b}
 \Delta u+2-2e^u-8(k-1)\|\beta_0+\db\eta\|^2e^{(k-1)u}=0,
\ene
and this is the first equation in the system \eqref{hit}, and furthermore:
\be\label{cp-2}
\int _X\Re\langle\beta_0+\db\eta,\overline{\partial}\ell\rangle e^{(k-1)u}=0, \text{\ \ \ for any $\ell \in A^0(E)$}.
\ene
By taking $\sqrt{-1}\ell$ instead of $\ell$ in \eqref{cp-2}, we find more generally that, at a {\cp} $(u,\eta)$ of $\Dcal$,
\be\label{cp-2b}
\int_X \langle\beta_0+\db\eta,\overline{\partial}\ell\rangle e^{(k-1)u}\ dA=0, \text{\ \ \ for any $\ell \in A^0(E)$}.
\ene
We need to show that \eqref{cp-2b} is equivalent to the condition that $e^{(k-1)u}\ast(\beta_0+\db\eta)$ is a {\hkd}, as stated by the second 
equation in the system \eqref{hit}. This is obtained as follows:
\bear
&&\int_X \langle\beta_0+\db\eta,\db\ell\rangle e^{(k-1)u}\ dA=0, \text{\ \ for any~} \ell\in A^0(E)\\
%&\Longleftrightarrow&\int_X \langle\beta_0+\db\eta,\db\ell\rangle e^{(k-1)u}\ dA=0, \text{\ \ for any~} \ell\in A^0(E)\\
&\Longleftrightarrow&\int_X \langle\db\ell,\beta_0+\db\eta\rangle e^{(k-1)u}\ dA=0, \text{\ \ for any~} \ell\in A^0(E)\\
&\Longleftrightarrow&\int_Xe^{(k-1)u} *(\beta_0+\db\eta)\wedge (\db\ell) dA=0, \text{\ \ for any~} \ell\in A^0(E)\\
&\Longleftrightarrow&\int_X \db[\ell e^{(k-1)u}*(\beta_0+\db\eta)]dA=\int_X \ell\db[e^{(k-1)u}*(\beta_0+\db\eta)]dA, \forall \ell\in A^0(E)\\
&\Longleftrightarrow&\int_X \ell\db[ e^{(k-1)u}*(\beta_0+\db\eta)]dA=0, \text{\ \ for any~} \ell\in A^0(E)\\
&\Longleftrightarrow&\db[e^{(k-1)u}*(\beta_0+\db\eta)]=0.
\end{eqnarray*}
Now we are left to show the regularity for the {\cp} $(u,\eta) \in \V$. We start with the following:
\vskip 0.1in
\underline{Claim}: If $\eta$ satisfies \eqref{cp-2b} with $u\in H^1(X)$, then $\eta\in W^{1,q}(X,E)$ for any $q\ge 1$.
\vskip 0.1in
\noindent
{\it Proof of the Claim}: To establish this, we consider the bundle $E\otimes K_X$ over X, which has a Hermitian inner product that arises from the Hermitian product defined on 
$X$. Since the surface is closed, we have that the space of holomorphic sections over it is a finite dimension vector space over $\C$ (see for instance Finiteness 
Theorems in \cite{Nar92}). Let us choose a basis $\{s_1, \cdots, s_N\}$ on this space of holomorphic sections such that
\beq
\int_X\langle\ast s_i, \ast s_j\rangle = \delta_i^j.
\eeq
Since \eqref{cp-2b} holds, by Weyl's regularity Lemma, we find that $e^{(k-1)u}\ast(\beta_0+\db\eta)$ is a holomorphic section (or a {\hkd} on $X$). Hence by 
the finiteness property of the space of holomorphic sections, for some $\alpha^i\in\C$, we have 
\beq
e^{(k-1)u}\ast(\beta_0+\db\eta) = \sum\limits_{i=1}^N\alpha^is_i.
\eeq
We denote the inverse map of the Hodge star $\ast$ by $\ast^{-1}$, then 
\be\label{dbeta1}
e^{(k-1)u}(\beta_0+\db\eta) = \sum\limits_{i=1}^N\alpha^i(\ast^{-1}s_i).
\ene
From this we have 
\be\label{dbeta}
\db\eta = -\beta_0 +e^{-(k-1)u}\sum\limits_{i=1}^N\alpha^i(\ast^{-1}s_i).
\ene
Since $e^u \in L^q(X)$ for all $q \ge 1$, and $\db$ is an elliptic operator, using trivialization and standard elliptic estimates, we deduce that 
$\eta \in W^{1,q}(X,E)$ for all $q \ge 1$. The Claim is proved.

Therefore we can apply this information in \eqref{cp-1b} and by elliptic regularity obtain that $u\in W^{2,q}(X)$ for all $q \ge 1$. In particular, we see that 
$u\in C^{1, b}(X)$ for some $b \in (0,1)$. Now elliptic regularity theory applied to equations \eqref{cp-1b} and \eqref{dbeta} combined with a bootstrapping 
argument allow us to obtain all the desired regularity for $(u,\eta)$.
\ep

%%%%%%%%%%%%%%%%%%%%%%%%%%%%%
\subsection{A Priori Estimates} Clearly {\df} is bounded from below by the value $4\pi(g-1)$. To analyze a {\msq}, we first provide some elementary estimates.
\bl\label{apriori}
For each $C> 0$, consider the sublevel set  
\beq
\Dcal^C:=\Big\{(u,\eta)\in\V: \Dcal(u,\eta) \le C \Big\}.
\eeq
Then we have:
\ben
\item
The set $\Big\{u: (u, \eta) \in \Dcal^C\Big\}$ is bounded in $H^1(X)$, and
\item
For any $a \in [1,2)$, there exists a constant $C_a > 0$ such that:
\be\label{beta-Lp}
\int_X\|\beta_0+\db\eta\|^a dA\le C_a, \ \ \forall (u,\eta)\in \Dcal^C.
\ene
\een
\el
\bp
We write $\beta = \beta_0+\db\eta$. Since the function $e^x - x$ is always positive, and by assumption we have $\Dcal(u,\eta) \le C$, 
then we find 
\be\label{gradbd}
\int_X \frac14|\nabla u|^2 dA \le C,
\ene
and
\be\label{gradbd2}
\int_X \|\beta\|^2e^{(k-1)u}dA \le \frac{C}{4},
\ene
and
\beq
\int_X \{e^{u} - u\}dA \le C.
\eeq
Furthermore we observe that for some positive constant $C_1$, we have:
\beq
e^{u} - u \ge |u| -C_1, \ \forall u\in \R.
\eeq
Thus, we conclude that,
\beq
\int_{X} |u| < C_2,
\eeq
for some suitable $C_2 > 0$. Hence by means of \eqref{gradbd} and Poincar\'e inequality we obtain the desired $H^1$-estimate.

To show \eqref{beta-Lp}, let $a \in [1,2), \ q = \frac{2}{a}$, and $b = \frac{(k-1)a}{2}$, we have
\bear
\int_X\|\beta\|^a dA&=& \int_X\|\beta\|^a e^{bu} e^{-bu}dA\\
&\le& \{\int_X\|\beta\|^{qa}e^{qbu}\}^{\frac1q}\{\int_X e^{-q'bu}\}^{\frac{1}{q'}}\\
&=& \{\int_X\|\beta\|^2 e^{(k-1)u}\}^{\frac{a}{2}}\{\int_X e^{-q'bu}\}^{\frac{2-a}{2}},
\eear
with $q' = \frac{2}{2-a}$.

For $(u,\eta)$ in the set $\Dcal^C$, we know that the first term on the right hand side is bounded by \eqref{gradbd2}, and also 
the second term is bounded by part (i) and Moser-Trudinger inequality.
\ep

%We note the following corollary to use it later:
%\bcor\label{corollary}
%For any sequence $(u_n, \eta_n)$ with $\Dcal(u_n, \eta_n)$ bounded, we have $\{u_n\}$ converges to some
%$u_0$ strongly in both $L^2(X)$ and $L^1(X)$. Moreover $\{e^{\pm u_n}\}$ converge strongly to $e^{\pm u_0}$
%in $L^q(X)$, for any $q \in [1,\infty)$. 
%\ecor

We prove the following important ``compactness" result:
\bpo\label{regu}
Let $(u_n,\eta_n) \in \V$ be a sequence satisfying 
\be\label{regu1}
\int_X e^{(k-1)u_n}\langle \beta_0+\db\eta_n,\db\ell\rangle =0, \text{\ \ for any~} \ell\in A^0(E).
\ene
Then if $u_n \rightharpoonup u$ weakly in $H^1(X)$, we have, along a subsequence, $\eta_n \to \eta$ in $W^{1,p}(X,E)$, with $\eta$ satisfying
\eqref{cp-2b}, namely, 
\beq
\int_X e^{(k-1)u}\langle \beta_0+\db\eta,\db\ell\rangle =0, \text{\ \ for any~} \ell\in A^0(E).
\eeq
\epo
\bp
Since $u_n  \xrightharpoonup{\, \, H^1(X)} u$, we know that $u_n$ is {\ubd} in $H^1(X)$, and {\upto} there holds, 
\beq
e^{(k-1)u_n} \xrightarrow{\, \, L^q} e^{(k-1)u}, \ \text{  for any } q\ge 1.
\eeq
By choosing $\ell = \eta_n$ 
in \eqref{regu1}, we have:
\bear
\int_X e^{(k-1)u_n}\|\db\eta_n\|^2\ dA& =& -\int_X e^{\frac{(k-1)u_n}{2}}\langle \beta_0, e^{\frac{(k-1)u_n}{2}}\db\eta_n\rangle\ dA\\
&\le& C_0\left(\int_X e^{(k-1)u_n} \right)^{\frac12}\left( \int_X e^{(k-1)u_n}\|\db\eta_n\|^2\right)^\frac12,
\eear
for some constant $C_0 > 0$. This implies 
\be\label{step1bd}
\int_X e^{(k-1)u_n}\|\db\eta_n\|^2\ dA \le C,
\ene
and so we have $\Dcal(u_n,\eta_n) \le C$ for some suitable $C>0$. Therefore  for some $a \in (1,2)$, by Lemma ~\ref{apriori}, we find a constant 
$C_a > 0$ such that $\int_X \|\db\eta_n\|^a \le C_a$. Hence by using Proposition ~\ref{PI}, we conclude that $\eta_n$ is bounded in $W^{1,a}(X,E)$, 
and therefore, {\upto}, $\eta _n  \rightharpoonup \eta$ in $W^{1,a}(X,E)$. In particular, for $a'=\frac{a}{a-1}$ the dual exponent of $a$, and for each 
fixed $\xi_0 \in L^{a'}(A^{0,1}(E))$, we have 
\be\label{xi}
\int_X\langle\db(\eta-\eta_n),\xi_0\rangle\ dA \to 0, \text{ as } n\to\infty,
\ene
since the map $\xi \longmapsto \int_X\langle\db\xi,\xi_0\rangle$ is a continuous linear map (by H\"older inequality).

Thus, for any smooth $\ell\in A^0(E)$, using \eqref{regu1}, we find:
\bear
\Big{|}\int_X e^{(k-1)u}\langle \beta_0+\db\eta,\db\ell\rangle\Big{|} &=& \Big{|}\int e^{(k-1)u}\langle \beta_0+\db\eta,\db\ell\rangle - 
e^{(k-1)u_n}\langle \beta_0+\db\eta_n,\db\ell\rangle\Big{|}\\
&\le& \Big{|}\int \big(e^{(k-1)u}-e^{(k-1)u_n}\big)\langle \beta_0+\db\eta_n,\db\ell\rangle\Big{|}\\
&& \ \ \ \ \ \ \ \ \ \ \ \ \ \ \ \ + \Big{|}\int e^{(k-1)u}\langle\db(\eta-\eta_n),\db\ell\rangle \Big{|}\\
&\le& \int \Big{|}e^{(k-1)u}-e^{(k-1)u_n}\Big{|}\|\beta_0+\db\eta_n\|\|\db\ell\|_\infty \\
&& \ \ \ \ \ \ \ \ \ + \Big{|}\int \langle\db(\eta-\eta_n),e^{(k-1)u}\db\ell\rangle \Big{|}\\
&\le&\|e^{(k-1)u}-e^{(k-1)u_n}\|_{L^{a'}}\|\beta_0+\db\eta_n\|_{L^a}\|\db\ell\|_\infty \\
&& \ \ \ \ \ \ \ \ \ +\Big{|}\int \langle\db(\eta-\eta_n),e^{(k-1)u}\db\ell\rangle \Big{|}\\
&\to& 0, \text{ as } n\to\infty.
\eear
Indeed, the first term in the last inequality goes to zero by the strong convergence of $e^{(k-1)u_n}$ to $e^{(k-1)u}$ in $L^{a'}$, and 
$\beta_0+\db\eta_n$ is bounded in $L^a$. The second term also goes to zero as a consequence of \eqref{xi} with $\xi_0 = e^{(k-1)u}\db\ell$. Now 
we have $(u, \eta) \in \V$ and \eqref{cp-2b} holds.

Now we are left to prove $\eta_n \to \eta$ in $W^{1,q}(X,E)$, for each $q>1$. To this end, we use \eqref{dbeta1} to write 
\be\label{alphai}
\alpha^i = \int_X\langle e^{(k-1)u}(\beta_0+\db\eta), (\ast^{-1}s_i)\rangle. 
\ene
and 
\be\label{alphani}
\alpha_n^i = \int_X\langle e^{(k-1)u_n}(\beta_0+\db\eta_n), (\ast^{-1}s_i)\rangle. 
\ene
Furthermore, for $1<a<2$, we have:
\bear
|\alpha^i-\alpha_n^i| & =& \Big{|}\int_X \langle e^{(k-1)u}(\beta_0+\db\eta)-e^{(k-1)u_n}(\beta_0+\db\eta_n), (\ast^{-1}s_i)\rangle \Big{|}\\
&\le&\|e^{(k-1)u}-e^{(k-1)u_n}\|_{L^{a'}}\|\beta_0+\db\eta_n\|_{L^a}\|(\ast^{-1}s_i)\|_\infty \\
&& \ \ \ \ \ \ \ \ \ +\Big{|}\int_X \langle\db(\eta-\eta_n),e^{(k-1)u}(\ast^{-1}s_i)\rangle \Big{|},
\eear
and as before we conclude that, $\alpha_n^i \to \alpha^i$, as $n \to \infty$. On the other hand,
\bear
\db(\eta_n-\eta) &=&  e^{-(k-1)u_n}\sum\limits_{i=1}^N\alpha_n^i(\ast^{-1}s_i) -  e^{-(k-1)u}\sum\limits_{i=1}^N\alpha^i(\ast^{-1}s_i)\\
&=& \big(e^{-(k-1)u_n}- e^{-(k-1)u}\big)\sum\limits_{i=1}^N\alpha_n^i(\ast^{-1}s_i)\\
&&\ \ \ \ \ \ \ \ \ \ \ \ +  e^{-(k-1)u}\sum\limits_{i=1}^N(\alpha_n^i-\alpha^i)(\ast^{-1}s_i),
\eear
and readily we derive that $\eta_n \xrightarrow {\, \, W^{1,q}} \eta$, for any $q \ge 1$. This concludes the proof.
\ep

%%%%%%%%%%%%%%%%%%%%
\subsection{The partial map}
\begin{definition}\label{Du}
For each fixed $u\in H^1(X)$, let us consider the following map: $D_u: W^{1,p}(X,E) \to \R$ with $\eta \mapsto \Dcal(u,\eta)$. 
 \end{definition}
This map will be very important in our strategy of proving the existence of a minimizer for {\df} in the space $\mathcal{V} = H^{1}(X) \times W^{1,p}(A^0(E))$, 
with $p > 2$. We first show this map has a unique minimum in $W^{1,p}(X,E)$.

\bt\label{du-map}
For each $u\in H^1(X)$, the map above $\Dcal_u$ admits a minimizer $\eta(u)$ which is its unique {\cp} in $W^{1,p}(X,E)$. Furthermore, $\eta(u)$ lies in 
$W^{1,q}(X,E)$, for all $q \ge 1$.
\et
\bp
Let $\eta_n \in W^{1,p}(X,E)$ be a minimizing sequence for the map $\Dcal_u$. By Lemma ~\ref{apriori}, for each $a \in (1,2)$, we have a constant 
$C_a > 0$ such that
\beq
\int_X\|\db\eta_n\|^a\ dA \le C_a.
\eeq
Therefore, for some $\eta \in W^{1,a}(X,E)$, we have $\eta_n \xrightharpoonup{\, \, W^{1,a}} {\eta}$. In addition, as in the proof of Proposition ~\ref{regu}, 
we find, for all smooth $\ell\in A^0(E)$,
%By Proposition ~\ref{regu}, applied with the 
%sequence $(u,\eta_n)$, we deduce that $\eta \in W^{1,p}(X,E)$, and $\eta_n \xrightarrow{\, \, W^{1,p}} {\eta}$
%We may choose the minimizing sequence $\\eta_n\in W^{1,p}(X,E)$ to satisfy:
%\beq \int_X e^{(k-1)u}\|\beta_0+\db\eta_n\|^2  \to \inf, \text{ as }\ n \to \infty, \eeq
% and for all smooth $\ell\in A^0(E)$,
% \beq \int_X e^{(k-1)u}\Re\langle\beta_0+\db\eta_n,\db\ell\rangle  \to 0, \text{ as }\ n \to \infty. \eeq
%We see that
%\bear |\int_X e^{(k-1)u}\Re\langle\beta_0+\db\eta_n,\db\ell\rangle &-& \int_X e^{(k-1)u}\Re\langle\beta_0+\db\eta,\db\ell\rangle| \\
%&\le&  \int_X e^{(k-1)u}\|\db(\eta_n-\eta)\|\|\db\ell\|\\ &\to& 0, \eear
%since $\eta_n \xrightharpoonup{\, \, W^{1,a}} {\eta}$. 
\beq
\int_X e^{(k-1)u}\langle\beta_0+\db\eta,\db\ell\rangle =0,
\eeq
and by regularity $\eta \in W^{1,q}(X,E)$, for all $q \ge 1$. Furthermore, 
from
\beq
\int_X e^{(k-1)u}\langle\db(\eta_n-\eta), \db(\eta_n-\eta)\rangle \ge 0,
\eeq
we get 
\bear
\int_X e^{(k-1)u}\langle\db\eta_n, \db\eta_n\rangle &\ge& \int_X \Big\{\langle\db(\eta_n-\eta), e^{(k-1)u}\db\eta\rangle+\langle e^{(k-1)u}\db\eta, \db(\eta_n-\eta)\rangle\Big\}\\
&& \ \ \ \ \ \ \ +\int_X e^{(k-1)u}\langle\db\eta, \db\eta\rangle. 
\eear
Since $\db(\eta_n-\eta) \xrightharpoonup{\, \, L^{a}} 0$, and $e^{(k-1)u}\db\eta \in L^{a'}$, so the first integral in the right hand side of the above inequality goes to 
zero as $n \to \infty$. Therefore we have 
\bear
\inf\limits_{\eta\in W^{1,p}(X,E)}\Dcal_u(\eta) &=&   \lim\limits_{n\to\infty}\int_X e^{(k-1)u}\|\beta_0+\db\eta_n\|^2\\
&\ge& \int_X e^{(k-1)u}\|\beta_0+\db\eta\|^2,
\eear
and therefore $\eta$ is a minimum for the map $\Dcal_u$.

We now observe that the map $\Dcal_u$ is strictly convex in $\eta$, and so $\eta$ is the unique {\cp} of $\Dcal_u$ in $W^{1,p}(X,E)$. 
\ep
Clearly, as in Theorem ~\ref{EL}, we see that, for each $u \in H^1(X)$, the unique {\cp} $\eta(u)$ of the map $\Dcal_u$ given by Theorem ~\ref{du-map} satisfies:
\beq
\int_X \langle e^{(k-1)u}(\beta_0+\db\eta(u)),\db\ell\rangle =0, \text{\ \ for any~} \ell\in A^0(E).
\eeq
Furthermore by the uniqueness, we deduce: 
\be\label{3.23}
\text{if }(u,\eta) \in \V, \text{and }\int_X \langle e^{(k-1)u}(\beta_0+\db\eta),\db\ell\rangle =0, \forall \ell\in A^0(E), \text{then }\eta=\eta(u). 
\ene
As a direct consequence of Proposition ~\ref{regu}, we deduce 
\bpo\label{A}
If $u_n \rightharpoonup u$ in $H^{1}(X)$, then $\eta(u_n) \rightarrow \eta(u)$ in $W^{1,q}(X,E)$, for all $q \ge 1$. 
\epo
We can now show:  
\bt\label{existence1}
The Donaldson functional $\Dcal$ admits a global minimum $(u_0, \eta_0)$ in $\V = H^{1}(X) \times W^{1,p}(X,E)$.
\et
\bp
Clearly the functional $\Dcal$ is bounded from below by the value $4\pi(g-1)$, and so we may consider a {\msq} $(u_n,\eta_n) \in \V$, such that 
$\Dcal(u_n,\eta_n)\to \inf\{\Dcal(u,\eta)\}$. By the definition of the map $\Dcal_u$ and Theorem ~\ref{du-map}, we have 
\beq
\Dcal (u_n, \eta_n) \ge \Dcal (u_n, \eta(u_n)),
\eeq
where $\eta(u_n)$ is the unique minimum for the map $\Dcal_{u_n}$. Therefore $(u_n, \eta(u_n))$ is also a {\msq} for {\df}. By part (i) of Lemma 
~\ref{apriori}, we can further assume that, $u_n \rightharpoonup u$ in $H^{1}(X)$, and therefore by Proposition ~\ref{A}, $\eta(u_n)$ converges to 
$\eta(u)$ in $W^{1,p}(X,E)$. Therefore, we find
\beq
\inf\Dcal =\lim\limits_{n\to\infty}\Dcal (u_n, \eta(u_n) \ge \Dcal (u, \eta(u)), 
\eeq
and so $(u,\eta(u))$ is a minimum for {\df}.
\ep

%%%%%%%%%%%%%%%%%%%%%%%%%%%%%%
%%%%%%%%%%%%%%%%%%%%%%%
\section{Second Variation of {\df}}

In this section, we study the second variation of {\df}. The main result in \S 4.1 is that, at a {\cp}, the Hessian is positive definite in the space 
$H^{1}(X)\times W^{1,2}(X,E)$ (Theorem ~\ref{hessian}). By additional estimates, we will prove in \S 4.2 that indeed every {\cp} is a strict 
local minimum for $\Dcal$ in the stronger space $\V$.

%%%%%%%%%%%%%%%%%%%%%%%%%%%%%%%%%%%%
\subsection{The Second Variation}
Recall from \eqref{df} that the functional $\Dcal(u,\eta)$ is defined as:
\beq
\Dcal(u,\eta)=\int_X\{\frac14|\nabla u|^2- u+e^u\}\ dA+4\int_X\|\beta_0+\db\eta\|^2e^{(k-1)u}\ dA,
\eeq
and at a {\cp} $(u,\eta) \in \V$ we have: 
\be\label{uk}
\Delta u+2-2e^u-8(k-1)\|\beta\|^2e^{(k-1)u}=0,
\ene
for $\beta = \beta_0+\db\eta$, and $\db\{e^{(k-1)u}\ast(\beta_0+\db\eta)\}=0$, that is
\be\label{eta-k}
\int_X \langle\beta_0+\db\eta, \db\ell\rangle e^{(k-1)u}\ dA =0, 
\ene
for all $\ell \in A^{0,1}(X,E)$.

We first observe the following consequence of the Proposition ~\ref{Bochner}:
\bl\label{Poi}
Let $E =\bigotimes\limits^{k-1} T^{1,0}_X$, and $\phi$ be a smooth function on $X$, then 
\be\label{Poi-k}
\int_X \|\db\ell\|^2e^{4(k-1)\phi}\ dA \ge (k-1)\int_X(\Delta\phi+1)\|\ell\|^2 e^{4(k-1)\phi}\ dA
\ene
holds for any $\ell \in W^{1,2}(X,E)$. 
\el 
\bp
We use the metric $g = e^{2\phi}g_X$ conformal to the {\hym} $g_X$. Then its Gaussian curvature can be calculated according to
\beq
K(g) = e^{-2\phi}(-\Delta\phi -1),
\eeq
where $\Delta$ is used {\wrt} the {\hym} $g_X$ which has constant sectional curvature $-1$.

Now using this Riemannian metric $g$, we have from \eqref{ineq1}
\beq
\int_X\langle\overline{\partial} \ell,\overline{\partial} \ell\rangle_g dA_g\geq -(k-1)\int_X K(g)\langle\ell,\ell\rangle_g dA_g.
\eeq
We change these inner products to be in terms of the {\hym} $g_X$. Since $E$ is the tensor product of $k-1$ copies of $T^{0,1}X$, and 
$\ell \in W^{1,2}(X,E)$, we have 
\beq
\langle\ell,\ell\rangle_g =e^{4(k-1)\phi}\langle\ell,\ell\rangle\
\eeq
and since $\db\ell$ is a $(0,1)$-form valued in $E$, we also have 
\beq
\langle\db\ell,\db\ell\rangle_g =e^{4(k-1)\phi}e^{-2\phi}\langle\db\ell,\db\ell\rangle.
\eeq
Finally we use the volume form $dA_g = e^{2\phi}dA$ to conclude the proof.
\ep
As a corollary, by taking $\phi = \frac{u}{4}$, we have 
\bcor\label{P3}
Let $E =\bigotimes\limits^{k-1} T^{1,0}_X$, and $(u,\eta)$ be a solution of the system of equations \eqref{uk} and \eqref{eta-k}, then 
\be\label{Poi-3}
\int_X \|\db\ell\|^2e^{(k-1)u} \ge 2(k-1)^2\int_X\|\beta\|^2\|\ell\|^2 e^{2(k-1)u}+ \frac{(k-1)}{2}\int_X \|\ell\|^2 e^{(k-1)u}
\ene
holds for any $\ell \in W^{1,2}(X,E)$.
\ecor
\bp
Since $u$ satisfies 
\beq
\Delta u+2-2e^u-8(k-1)\|\beta\|^2e^{(k-1)u}=0,
\eeq
we choose $\phi = \frac{u}{4}$ in Lemma ~\ref{Poi} to find:
\bear
\int_X \|\db\ell\|^2e^{(k-1)u} &\ge& \frac{(k-1)}{2}\int_X\Big\{e^u+1+4(k-1)\|\beta\|^2e^{(k-1)u}\Big\}\|\ell\|^2 e^{(k-1)u}\\
&\ge&2(k-1)^2\int_X\|\beta\|^2\|\ell\|^2 e^{2(k-1)u} + \frac{(k-1)}{2}\int_X \|\ell\|^2 e^{(k-1)u}.
\eear
\ep
Our main result in this subsection is the following:
\bt\label{hessian}
At any {\cp} $(u,\eta)$, setting $\beta = \beta_0+\db\eta$, the second variation 
of {\df} $\Dcal(u,\eta)$ is strictly positive definite. More specifically we have, for all $v \in H^1(X)$ and all $0\not=\ell\in W^{1,p}(X,E)$,
\begin{eqnarray}\label{hes1}
\Dcal''_{(u,\eta)}(v,\ell) &=&  \int_X\Big\{e^uv^2 + 4\|(k-1)v\beta+\db\ell\|^2e^{(k-1)u}\Big\}dA \notag\\
&&\ \ \ \ \ \ + 2\int_X\Big\{\Big\|2(k-1)e^{(k-1)u}\|\ell\|\beta - \frac{\db v\otimes\ell}{\|\ell\|}\Big\|^2\big\}dA +\mathcal{R}_{(u,\eta)}(v,\ell),
\end{eqnarray}
where
\begin{eqnarray}\label{R}
\mathcal{R}_{(u,\eta)}(v,\ell) &=& 4\int_X\big\{\|\db\ell\|^2\ e^{(k-1)u}- 2(k-1)^2\|\beta\|^2e^{2(k-1)u}\|\ell\|^2\big\}dA \\
 &\ge& 2(k-1)\int_X e^{(k-1)u}\|\ell\|^2 dA. \notag
\end{eqnarray}
%In particular, we have:\be\label{hes2}\Dcal''_{(u,\eta)}(v,\ell) \ge \int_X\Big\{e^uv^2 + 2(k-1)e^{(k-1)u}\|\ell\|^2\Big\}\ dA,\ene
\et
\bp
At a {\cp} $(u,\eta)$, we have for any tangent vector $(v,\ell)$, with $v \in H^1(X)$ and $\ell \in W^{1,p}(X,E)$, 
the second variation can be computed as follows:
\bear
\Dcal''(v,\ell) &&= \int_X\Big\{\frac{|\nabla v|^2}{2}+e^uv^2\Big\}\\ 
&&\ \ \ \ \ +4\int_X \Big\{(k-1)^2\|\beta\|^2 v^2+4(k-1)\Re\langle\beta,\overline{\partial}\ell\rangle v+2\|\db\ell\|^2\Big\} e^{(k-1)u}\\
&&= A + B +  \int_X\Big\{e^uv^2\Big\},
\eear
where we write 
\beq
A = 4\int_X \Big\{(k-1)^2\|\beta\|^2 v^2+2(k-1)\Re\langle\beta,\overline{\partial}\ell\rangle v+\|\overline{\partial}\ell\|^2\Big\} e^{(k-1)u},
\eeq
and 
\beq
B = \int_X\Big\{\frac{|\nabla v|^2}{2}\Big\}+4\int_X \Big\{2(k-1)\Re\langle\beta,\overline{\partial}\ell\rangle v+\|\overline{\partial}\ell\|^2\Big\} e^{(k-1)u}.
\eeq
Clearly we have: 
\begin{eqnarray}\label{Aa}
A &=& 4\int_X \langle (k-1)v\beta+\db\ell,(k-1)v\beta+\db\ell\rangle e^{(k-1)u}\\
&\ge& 0 \nonumber.
\end{eqnarray}

Now let us consider the term $B$. Since $(u,\eta)$ is a {\cp}, we have:
\beq
\int_X \langle\beta,\db(v\ell)\rangle e^{(k-1)u} = 0. 
\eeq
%To see this 
%\begin{eqnarray*}
%\int_X Re\langle\beta,\db\ell\rangle ve^{(k-1)u}&=&\int_X Re\langle\db\ell,\beta\rangle ve^{(k-1)u} \\
%&=&\int_X Re(\overline{\partial \ell\wedge *\beta}ve^{(k-1)u})\\
%&=&\int_X Re(\db[\ell v*\beta])e^{(k-1)u}-\int_X Re(\db(v*\beta) \ell)e^{(k-1)u} \\
%&=&-\int_X Re(\db(v*\beta) \ell), \text{\ \ since\ \ } \db*\beta=0\\
%&=&-\int_X Re(\overline{\partial}(v)\wedge*\beta \ell)\\
%&=&-\int_X Re\langle\beta,\overline{\partial}(v) \ell\rangle.
%\end{eqnarray*}
Then by Leibniz's rule for the operator $\db_E=\db$ (see for instance \cite{Wel08}), we have 
\beq
\db_E(v\ell) = (\db v)\otimes\ell + v\db_E\ell,
\eeq
where $\db v = \partial_{\bar{z}}(v)d\bar{z} \in A^{0,1}(X)$ for a function $v$. Therefore 
\beq
\int_X \langle\beta, v(\db\ell)\rangle e^{(k-1)u} = -\int_X \langle\beta,(\db v)\otimes\ell\rangle e^{(k-1)u}.
\eeq
Noting that 
\beq
\|\db v\|^2 = \frac14|\partial_xv+\sqrt{-1}\partial_yv|^2\|d\bar{z}\|^2 = \frac14|\nabla v|^2,
\eeq 
we can express $B$ as follows:
\beq
B = 2\int_X\Big\{\|\db v\|^2 -4(k-1)\Re\langle\beta,\db v\otimes\ell\rangle e^{(k-1)u}+2\|\db\ell\|^2\ e^{(k-1)u}\Big\}.
\eeq
Since $\|\db v\|^2 = \|\frac{\db v\otimes\ell}{\|\ell\|}\|^2$, the above equality is equivalent to 
\be\label{newB}
B = 2\int_X\Big\|\frac{\db v\otimes\ell}{\|\ell\|}-2(k-1)e^{(k-1)u}\|\ell\|\beta\Big\|^2\ dA+\mathcal{R}_{(u,\eta)}(v,\ell),
\ene
with the remainder term 
\beq
\mathcal{R}_{(u,\eta)}(v,\ell) = 4\int_X\Big\{\|\db\ell\|^2\ e^{(k-1)u}- 2(k-1)^2\|\beta\|^2e^{2(k-1)u}\|\ell\|^2\Big\}dA .
 %&\ge& 2(k-1)\int_X e^{(k-1)u}\|\ell\|^2 dA. 
\eeq
By \eqref{Poi-3}, we have established 
\beq
\int_X \|\db\ell\|^2\ e^{(k-1)u} \ge 2(k-1)^2\int_X\|\beta\|^2e^{2(k-1)u}\|\ell\|^2 + \frac{(k-1)}{2}\int_X e^{(k-1)u}\|\ell\|^2.
\eeq
Therefore we have 
\beq
\mathcal{R}_{(u,\eta)}(v,\ell) \ge 2(k-1)\int_X e^{(k-1)u}\|\ell\|^2 dA.
\eeq
This completes the proof. %Alternatively, using Cauchy-Schwarz inequality, we have 
%\be\label{Bb}B \ge 2\int_X\big\{\|\db v\| - 2(k-1)\|\beta\|e^{(k-1)u}\|\ell\|\big\}^2 + 2(k-1)\int_X e^{(k-1)u}\|\ell\|^2,\ene
%and the conclusion follows.
\ep
%%%%%%%%%%%%%%%%%%%%%%%%%%%%%%%%
\subsection{Every {\cp} is a strict local minimum}
Unfortunately the estimates in Theorem ~\ref{hessian} do not directly imply that every {\cp} is a strict local minimum of $\Dcal$ 
in the stronger norm of $\V$. As we shall see in Proposition ~\ref{ast}, these estimates ensures only that this is true in the weaker 
(but natural) space $H^1(X)\times W^{1,2}(X,E)$.
\bl\label{H1norm}
Let $\tilde\beta$ be a given continuous $(0,1)$-form valued in $E$. Then there exists $\sigma = \sigma(\tilde\beta)>0$ such that 
$\forall\ (f,F)\in L^2(X)\times L^2(A^{0,1}(X,E))$, 
\be\label{H1H1}
\int_X f^2\ dA + \int_X \|f\tilde\beta+F\|^2\ dA \ge \sigma\int_X (f^2+ \|F\|^2)\ dA.
\ene
\el
\bp
It is equivalent to prove 
\be\label{equiv}
\inf\limits_{\|f\|_{L^2}^2+\|F\|_{L^2}^2 =1}\Big\{\int_X f^2 + \|f\beta+F\|^2\ dA\Big\}> 0.
\ene
We argue by contradiction. Let $(f_n,F_n)\in L^2(X)\times L^2(A^{0,1}(X,E))$ be a sequence with $\|f_n\|_{L^2}^2+\|F_n\|_{L^2}^2 =1$ such that 
$\int_X (f_n^2+ \|f_n\tilde\beta+ F_n\|^2)\ dA\to 0$, as $n\to\infty$. Thus we have $\int_X f_n^2\ dA\to 0$ and 
$\int_X \|f_n\tilde\beta+ F_n\|^2\ dA\to 0$. Therefore,
\beq
\|F_n\|^2 = \|F_n+f_n\tilde\beta-f_n\tilde\beta\|^2 \le \|F_n+f_n\tilde\beta\|^2 + \|f_n\tilde\beta\|^2 \to 0,
\eeq
and this contradicts the assumption: $\|f_n\|_{L^2}^2+\|F_n\|_{L^2}^2 =1$.
\ep
As a consequence, we also obtain:
\bpo\label{ast}
For every {\cp} $(u,\eta)$ of $\Dcal$, there exists $\sigma = \sigma(u,\eta) > 0$ such that, $\forall\ (v, \ell) \in H^1(X)\times W^{1,2}(X,E)$,
\beq
\Dcal''_{(u,\eta)}(v,\ell) \ge \sigma\big( \|v\|_{H^1(X)}^2+\|\ell\|_{W^{1,2}(A^{0}(E))}^2 \big).
\eeq
\epo
\bp
The second variation $\Dcal''_{(u,\eta)}(v,\ell)$ at a {\cp} $(u,\eta)$ is given in the expression \eqref{hes1}. Since $(u, \eta)$ is smooth, we apply 
Lemma ~\ref{H1norm} with $\tilde\beta=(k-1)\beta=(k-1)(\beta_0+\db\eta)$ fixed, and $(f,F) =(v,\db\ell) \in L^2(X)\times L^2(A^{0,1}(X,E))$, to find
\bear
\int_X\big\{e^uv^2+4\|(k-1)v\beta+\db\ell\|^2e^{(k-1)u}\big\}dA&\ge&\sigma_1 \int_X\big\{v^2+\|(k-1)v\beta+\db\ell\|^2\big\}\ dA\\
&\ge& \sigma_1 \int_X\big\{v^2+\|\db\ell\|^2\big\}\ dA.
\eear
Now we inspect the last two terms in \eqref{hes1}. From \eqref{R}, and the fact that $e^{(k-1)u}$ is smooth, we have 
\beq
\mathcal{R}_{(u,\eta)}(v,\ell) \ge C\int_X \|\ell\|^2 dA.
\eeq
 Next we apply Lemma ~\ref{H1norm} again, with $\tilde\beta=2(k-1)e^{(k-1)u}\beta$ and $(f,F) =(\|\ell\|,\frac{\db v\otimes\ell}{\|\ell\|})$, to get 
\bear
2\int_X\Big\{\Big\|2(k-1)e^{(k-1)u}\|\ell\|\beta - \frac{\db v\otimes\ell}{\|\ell\|}\Big\|^2\big\}\ dA +\mathcal{R}_{(u,\eta)}(v,\ell) \ge \sigma_2 \int_X(\|\db v\|^2+\|\ell\|^2).
\eear
The proof is complete.
\ep
Clearly Proposition ~\ref{ast} does not suffice to show that a {\cp} of $\Dcal$ is a strict local minimum in $\V$. To this purpose, we establish the following estimate:
\bl\label{lemma2}
Let $u_0 \in H^1(X), \ \eta_0 := \eta(u_0)$, and $b> 0$. Then there exists a constant $C_0 = C_0(b,\|u_0\|_{H^1(X)}) > 0$ such that: 
\be\label{etabound}
\|\db\eta(u)-\db\eta_0\|^2_{L^p} \le C_0\big\{\|\db\eta(u)-\db\eta_0\|^2_{L^2}+ \|u-u_0\|^2_{L^2(X)} \big\}
\ene
holds for all $u \in H^1(X)$ with $\|u\|_{H^1(X)} < b$.
\el
\bp
Recall that the space 
\beq
\Cal_k(X) = \{\alpha \in A^{1,0}(X,E^*): \db\alpha = 0\} 
\eeq
is the space of {\hkd}s on $X$, and it is a finite dimensional vector space. Therefore all norms on $\Cal_k(X)$ are equivalent. Hence for any $q \ge 1$, there is a 
constant $C_q > 0$ such that 
\be\label{Aq}
\|\alpha\|_{L^q} \le C_q\|\alpha\|_{L^1}, \ \ \forall \alpha \in \Cal_k(X),
\ene
and we note $e^{(k-1)u}*(\beta_0+\db\eta(u)) \in \Cal_k(X)$. Furthermore, by \eqref{MTI2}, for each $q \ge 1$, there is a 
suitable constant $B_q > 0$, depending only on $\|u_0\|_{H^1(X)}$, $q$ and $b$, such that 
\be\label{Bq}
\|e^{(k-1)u} - e^{(k-1)u_0}\|_{L^q} \le B_q\|u-u_0\|_{L^2(X)}.
\ene
We use
\bear
\db\eta(u)-\db\eta_0 &=& e^{-(k-1)u}\Big\{e^{(k-1)u}(\beta_0+\db\eta(u)) - e^{(k-1)u_0}(\beta_0+\db\eta_0)\Big\} \\
&&\ \ \ \ \ \ \ \ +(e^{-(k-1)u}-e^{-(k-1)u_0})e^{(k-1)u_0}(\beta_0+\db\eta_0)
\eear
to estimate: 
\bear
\|\db\eta(u)-\db\eta_0\|_{L^p} &\le& \|e^{-(k-1)u}\{e^{(k-1)u}(\beta_0+\db\eta(u)) - e^{(k-1)u_0}(\beta_0+\db\eta_0)\}\|_{L^p} \\
&&\ \ \ \ \ \ \ \ +\|(e^{-(k-1)u}-e^{-(k-1)u_0})e^{(k-1)u_0}(\beta_0+\db\eta_0)\|_{L^p}.
\eear
Using H\"{o}lder inequality and the fact that the operator $\ast: A^{0,1}(X,E) \to A^{1,0}(X,E^*)$ is an isometry, and setting 
$\alpha_0 = \ast e^{(k-1)u_0}(\beta_0+\db\eta_0)$ and $\alpha = \ast e^{(k-1)u}(\beta_0+\db\eta(u))$, we have:  
\bear
\|\db\eta(u)-\db\eta_0\|_{L^p} &\le& \|e^{-u}\|_{L^{2(k-1)p}}\|\alpha - \alpha_0\|_{L^{2p}} \\
&&\ \ \ \ \ \ \ \ +\|(e^{-(k-1)u}-e^{-(k-1)u_0})e^{(k-1)u_0}(\beta_0+\db\eta_0)\|_{L^p}.
\eear
Therefore, by means of \eqref{MTI1}, the estimates \eqref{Aq} and \eqref{Bq}, there exists a constant $C_1>0$ (depending only on 
$b, p, \|u_0\|_{H^1}$) such that 
\be\label{Cp}
\|\db\eta(u)-\db\eta_0\|_{L^p} \le C_1\big(\|u-u_0\|_{L^2}+\|\alpha - \alpha_0\|_{L^1}\big).
\ene
Writing
\beq
\alpha-\alpha_0 = \ast\Big\{e^{(k-1)u}(\db\eta(u)-\db\eta_0) +(e^{(k-1)u}-e^{(k-1)u_0})(\beta_0+\db\eta_0)\Big\},
\eeq
we deduce 
\be\label{Dp}
\|\alpha - \alpha_0\|_{L^1} \le \|e^{(k-1)u}\|_{L^2}\|\db\eta(u)-\db\eta_0\|_{L^2} + C_2\|u-u_0\|_{L^2},
\ene
where, again, $C_2 = C_2(b, p, \|u_0\|_{H^1}) > 0$. Thus, from \eqref{Cp} and \eqref{Dp}, we obtain a constant 
$C_0 = C_0(b, p, \|u_0\|_{H^1}) > 0$, such that 
\beq
\|\db\eta(u)-\db\eta_0\|_{L^p} \le C_0\big\{\|\db\eta(u)-\db\eta_0\|_{L^2}+ \|u-u_0\|_{L^2(X)} \big\}
\eeq
for all $u \in H^1(X)$ with $\|u\|_{H^1(X)} < b$, and \eqref{etabound} is established.
\ep
Now we are ready to prove that each {\cp} of $\Dcal$ in $\V$ is a strict local minimum. To start, we show the following:
\bt\label{Lem3}
Let $p_0 =(u_0,\eta_0)$ be a {\cp} for $\Dcal$. Then there exist $\delta_0>0$ and $t_0>0$ such that 
$\Dcal(u,\eta(u)) \ge \Dcal(u_0,\eta_0) + t_0\|(u,\eta(u))-p_0\|^2_{\V}$ holds whenever $\|u-u_0\|_{H^1(X)} < \delta_0$.
\et
\bp
Given $r > 0$, let $B_r(p_0)$ in $(\V, \|\cdot\|_{\V})$ be the ball centered at $p_0$ of radius $r$. By the continuity of the map $u\to (u,\eta(u))$ 
from $H^1(X)$ to $\V$,  there exists some $\delta_r >0$, such that $(u,\eta(u)) \in B_{r}(p_0)$, whenever $\|u-u_0\|_{H^1} < \delta_r$. In particular, 
$0 < \delta_r < r$.

We apply Taylor expansion in the fixed ball $B_r(p_0)$, and for $(u,\eta) \in B_{r}(p_0)$, in virtue of Proposition ~\ref{ast} (with 
$\sigma_0 = \sigma(u_0,\eta_0)$), we can write 
\bear
\Dcal(u,\eta) &=& \Dcal(p_0) + \frac{\Dcal''_{p_0}[u-u_0,\eta-\eta_0]}{2} +o(\|(u, \eta) -p_0\|_{\V}^2)\\
 &\ge& \Dcal(p_0) + \frac{\sigma_0}{2}(\|u-u_0\|_{H^1}^2+\|\eta-\eta_0\|_{W^{1,2}}^2) +o(\|(u, \eta) -p_0\|_{\V}^2), \text{as } r\to 0^{+}.
\eear
On the other hand, when $\eta = \eta(u)$ then by Lemma ~\ref{lemma2} and Poincar\'{e} inequality, we obtain
\bear
\|u-u_0\|_{H^1}^2+\|\eta(u)-\eta_0\|_{W^{1,2}}^2 &\ge& \frac12\|u-u_0\|_{H^1}^2+\frac{1}{2C_0}\|\db\eta(u)-\db\eta_0\|^2_{L^p} \\
&\ge& \alpha_p\big(\|u-u_0\|_{H^1}^2+ \|\eta(u)-\eta_0\|_{W^{1,p}}^2\big),
\eear
with a suitable $\alpha_p > 0$ for any $u \in H^1(X): \|u-u_0\|_{H^1} < \delta_r$, and for any $r>0$ sufficiently small. As a consequence, we have:
\beq
\Dcal(u,\eta(u)) \geq \Dcal(p_0) + \Big(\frac{\sigma_0\alpha_p}{2}+o(1)\Big)\|(u,\eta(u))-p_0\|_{\V}^2, \text{ as } r\to 0^{+}. 
\eeq
Thus, by choosing $t_0 = \frac{\sigma_0\alpha_p}{4}$, we find $r_0 > 0$ sufficiently small and corresponding $\delta_0=\delta_{r_0} > 0$, such that, 
for $\|u-u_0\|_{H^1} < \delta_0$, we have:
\beq
\Dcal(u,\eta(u)) \ge \Dcal(u_0,\eta_0) + t_0\|(u,\eta(u))-p_0\|_{\V}^2,
\eeq
as claimed.
\ep
Consequently, we deduce:
\bcor\label{sphere}
Let $p_0=(u_0,\eta_0)$ be a {\cp} of {\df} $\Dcal$. Then $p_0$ is a strict local minimum for $\Dcal$ in $\V$. More precisely, for suitable 
$\delta_0>0$ sufficiently small it holds: 
\be\label{sphere2}
\Dcal(u,\eta) > \Dcal(u_0,\eta_0), \ \forall (u,\eta) \in \V \text{ with } \|u-u_0\|_{H^1} < \delta_0, \text{ and } (u,\eta) \not= (u_0,\eta_0).
\ene
\ecor
\bp
For any given {\cp} $p_0 =(u_0,\eta_0)$ of $\Dcal$, let $\delta_0> 0$ and $t_0>0$ be as given in Theorem ~\ref{Lem3}, so that 
\be\label{sphere1b}
\Dcal(u,\eta) \ge \Dcal(u,\eta(u)) \ge \Dcal(u_0,\eta_0) + t_0\|(u,\eta(u))-p_0\|_{\V}^2
\ene
holds, whenever $\|u-u_0\|_{H^1} < \delta_0$.

Now we assume that $(u,\eta)\not= p_0$, then there are two cases to consider. In the case $u \not= u_0$, then from \eqref{sphere2}, we 
deduce that $\Dcal(u,\eta) > \Dcal(u_0,\eta_0)$, as claimed. In the other case where $u =u_0$, then necessarily we have 
$\eta\not=\eta_0=\eta(u_0)$. Since $\eta_0$ is the unique strict minimum of the partial map $\Dcal_{u_0}$, in this case we have 
\beq
\Dcal(u,\eta) = \Dcal(u_0,\eta) =\Dcal_{u_0}(\eta) > \Dcal_{u_0}(\eta_0)= \Dcal(u_0,\eta_0).
\eeq
In conclusion \eqref{sphere2} holds and the proof is completed.
\ep
%%%%%%%%%%%%%%%%%%%%%%%%%%%%%%%%%%%%%%%%%%%%%%
\section{A ``Weaker" {\PS} Condition, Ekeland Principle and Uniqueness of the Critical Point} 

As pointed out in the introduction, a functional may have many {\cp}s which are all strict local minima. Thus to prove that $\Dcal$ admits a unique 
{\cp} (i.e. its global minimum), our approach is to assume (by contradiction) the existence of more strict local minima for $\Dcal$ then arrive to a 
contradiction by a ``{\mpa}" construction (\cite{AR73}) that yields to an additional {\cp} of $\Dcal$, which however is \underline{not} a local minimum. 
But to successfully carry out such program, we need {\df} to satisfy the following (well known) {\PSc}: if a sequence $(u_n,\eta_n) \in \V$ satisfies:
\be\label{PSq}
\Dcal(u_n,\eta_n) \to c, \text{ and } \|\Dcal'_{(u_n,\eta_n)}\|\to 0, \ \text{ as } n\to\infty,
\ene
then up to a subsequence, $(u_n,\eta_n)$ (called a {\PS} sequence) converges strongly in $\V$. Unfortunately, it is not obvious to check such 
property, since \eqref{PSq} does not provide any reasonable control of the component $\eta_n$ in the space $W^{1,p}(X,E)$ when $p > 2$. In fact 
as seen above, the best we can hope for is a uniform estimate in $W^{1,a}(X,E)$, with $1< a < 2$.

On the other hand, if we assume a priori that $\{\eta_n\}$ is bounded in $W^{1,p}(X,E)$, a ``weaker" form of the {\PSc} holds (see 
Lemma ~\ref{almostPS}), which suffices for our purpose. Indeed, by means of the Ekeland Principle (Theorem ~\ref{Eke} below) we are able to obtain an 
``ad-hoc" {\PS} sequence satisfying such additional uniform bound for the component $\eta_n$. 

%%%%%%%%%%%%%%%%%%%%%%%%
\subsection{A ``weaker" {\PSc}}
We start with the following lemma:
\bl\label{almostPS}
Let $(u_n,\eta_n) \in \V$ be a {\PS} sequence satisfying \eqref{PSq}. If $\{\eta_n\}$ is uniformly bounded in $W^{1,p}(X,E)$, then {\upto}, as 
$n\to\infty$, we have
\ben
\item 
$u_n \to u$ strongly in $H^1(X)$, and $\eta_n\to \eta$ strongly in $W^{1,2}(X,E)$;
\item
$\Dcal(u_n,\eta_n) \to \Dcal(u,\eta)$ and $\Dcal'_{(u,\eta)}= 0$.
\een
In particular, $c$ is a critical value for $\Dcal$ with corresponding {\cp} $(u,\eta)$.
\el
\bp
Since by assumption we have $\Dcal (u_n, \eta_n) \le C$ and $\|\eta_n\|_{W^{1,p}(X,E)}\le C$, then, {\upto}, we can assume: 
\beq
u_n \xrightharpoonup{\, \, H^1} u,\ \ \eta_n\xrightharpoonup{\, \, W^{1,p}(X,E)}\eta;
\eeq
and in addition:  
\be\label{5.2}
u_n \xrightarrow{L^\alpha} u, \ \ \ \ \text{and}\ \ \ e^{u_n} \xrightarrow{L^\alpha} e^{u}, \ \ \forall\alpha\ge1,
\ene
as $n\to\infty$.

Furthermore, by assumption we have: 
\be\label{psc1}
|\Dcal'_{(u_n,\eta_n)}(\xi, \ell)| = o(1)\|(\xi, \ell)\|_{\V} \to 0, \text{ as } n\to\infty, \ \ \forall (\xi,\ell) \in \V.
\ene
Thus, by arguing as in Proposition ~\ref{regu}, we immediately derive that $\Dcal'_{(u,\eta)}= 0$, and so $(u,\eta)$ is a (\underline{smooth}) {\cp} 
of $\Dcal$. Moreover, by using this information together with \eqref{psc1} with $\xi =0$ and 
$\ell = \eta_n-\eta$ (uniformly bounded in $W^{1,p}(X,E)$), we obtain:
\beq
\int_X e^{(k-1)u_n}\|\db(\eta_n-\eta)\|^2 dA = o(1), \text{ as } n \to \infty,
\eeq 
and consequently,
\beq
\int_X e^{(k-1)u}\|\db(\eta_n-\eta)\|^2 dA = o(1), \text{ as } n \to \infty.
\eeq 
Therefore, (by Poincar\'e inequality), $\eta_n\to \eta$ strongly in $W^{1,2}(X,E)$, and 
\beq
\int_X e^{(k-1)u_n}\|\beta_0+\db\eta_n\|^2 dA \to \int_X e^{(k-1)u}\|\beta_0+\db\eta\|^2 dA, \text{ as } n \to \infty.
\eeq 
Next we choose $\xi_n =u_n-u$ and $\ell = 0$ in \eqref{psc1} to find 
\beq
|\Dcal'_{(u_n,\eta_n)}(u_n-u,0)|= o(1), \text{ as }n \to \infty.
\eeq
This means,
\beq
\int_X \frac{\nabla u_n}{2}\nabla(u_n-u) + (e^{u_n}-1)(u_n-u)+4(k-1)\|\beta_0+\db\eta_n\|^2e^{(k-1)u_n}(u_n-u) \to 0, 
\eeq 
as $n\to +\infty$. Thus, we have 
\bear
\int_X \frac12|\nabla(u_n-u)|^2 &=& \int_X(u_n-u) - \int_X e^{u_n}(u_n-u) \\
&&\ \ \ \ \ \ \ \ -4(k-1)\int_X\|\beta_0+\db\eta_n\|^2e^{(k-1)u_n}(u_n-u) + o(1),
\eear
and the right hand side goes to $0$ as $n\to+\infty$. Since by \eqref{5.2} we know that 
\beq
\int_X e^{u_n}(u_n-u) \to 0, \text{ in }L^q(X), \forall q \ge 1,
\eeq
and so by using H\"older inequality, and by recalling that $\eta_n$ is uniformly bounded in $W^{1,p}(A^0(E))$,
\bear
\int_X\|\beta_0+\db\eta_n\|^2e^{(k-1)u_n}(u_n-u) &\le& \left(\int_X\|\beta_0+\db\eta_n\|^{p}\right)^\frac2p
\|e^{(k-1)u_n}(u_n-u)\|_{\frac{p-2}{p}}\\
&\to& 0, \text{ as }\  n\to +\infty.
\eear
Hence, $u_n\to u$ strongly in $H^1(X)$. Consequently, $\Dcal(u_n,\eta_n) \to \Dcal(u,\eta) = c$, and the proof is complete.
\ep
\br
We note that, even under the stronger assumption of Lemma ~\ref{almostPS}, we do not know whether or not the sequence $(u_n, \eta_n)$ 
converges in $\V$.
\er
%%%%%%%%%%%%%%%%%%%%%%%
\subsection{The Ekeland principle}
Let us first recall the Ekeland's $\epsilon$-variational principle as follows:
\bt(\cite{AE84})\label{Eke}
Let $(Y, d)$ be a complete metric space, and $F: Y \to \R$ a nonnegative and lower semi-continuous functional. Let there 
be given $\epsilon > 0$, and $\gamma_\epsilon^0 \in Y$, such that 
%\be\label{eke0}
$F(\gamma_\epsilon^0) \le \epsilon + \inf F$.
%\ene
Then there is some point $\gamma_\epsilon \in Y$ such that 
\be\label{eke1}
F(\gamma_\epsilon) \le F(\gamma_\epsilon^0),
\ene
\be\label{eke2}
d(\gamma_\epsilon,\gamma_\epsilon^0) \le \sqrt{\epsilon},
\ene
and 
\be\label{eke3}
F(\gamma) \ge F(\gamma_\epsilon) - \sqrt{\epsilon}d(\gamma,\gamma_\epsilon), \text{ for all } \gamma \in Y.
\ene
\et
We will use the Ekeland principle to prove an important lemma in a more general form than we need but may be of independent 
interest. To this end, we consider two distinct points $P_1$ and $P_2$ in some Banach space $(V, \|\cdot\|)$, and set
\be\label{path0}
\Pat = \big\{\gamma \in C^0([0,1], V): \gamma(0) = P_1,\text{ and } \gamma(1) =P_2 \big\}.
\ene
Clearly $\Pat$ is not empty, as it contains the path: $\gamma(t) = (1-t)P_1+tP_2$, and $(\Pat, d)$ is a complete 
metric space equipped with the metric $d(\gamma_1,\gamma_2) = \max\limits_{t \in [0,1]}\|\gamma_1(t) - \gamma_2(t)\|$.
\bl\label{grad-D}
Let $(V, \|\cdot\|)$ be a Banach space which is uniformly convex and $J: V\to \R$ be a $C^1$-function on $V$. Suppose there exist 
$\epsilon > 0$ and $\gamma_\epsilon \in \Pat$ such that
\be\label{eke0b}
\max\limits_{t \in [0,1]} J(\gamma(t)) \ge \max\limits_{t \in [0,1]} J(\gamma_\epsilon(t))- \sqrt{\epsilon}d(\gamma,\gamma_\epsilon)
\ene
holds for all $\gamma \in \Pat$. Set 
\be\label{2.7}
T_\epsilon = \big\{\tilde{t} \in [0,1]: J(\gamma_\epsilon(\tilde{t})) =  \max\limits_{t \in [0,1]} J(\gamma_\epsilon(t)) \big\}.
\ene
If $T_\epsilon \subset\subset (0,1)$, i.e., compactly contained in $(0,1)$, then there is $t_\epsilon \in T_\epsilon$ such that 
\be\label{2.8}
\|J'_{\gamma_\epsilon(t_\epsilon)}\|_{\ast} \le \sqrt{\epsilon}.
\ene
\el
\bp
We define $F(\gamma) = \max\limits_{t \in [0,1]} J(\gamma(t))$. Let $\rho_\epsilon: [0,1]\to [0,1]$ be a continuous function such that 
$\rho_\epsilon(0) = \rho_\epsilon(1) =0$, and $\rho_\epsilon(t) \equiv 1$, $\forall\ t \in T_\epsilon$. This cut-off function exists since 
$T_\epsilon \subset\subset (0,1)$.

Since $V$ is uniformly convex, it is reflexive by a theorem of Milman-Pettis (\cite{Bre83}), and its bi-dual $V^{\ast\ast}$ is also uniformly 
convex (since the canonical map $V\to V^{\ast\ast}$ is an isometry). In particular, given $f\in V^\ast$, there exists a unique 
$\hat{f} \in V^{\ast\ast}$ satisfying $\hat{f}(f) = \|f\|_\ast^2$ and $\|\hat{f}\| = \|f\|_\ast$, where $\|\cdot\|_\ast$ is the norm on $V^\ast$. This 
gives a well-defined ``duality map" $V^\ast\to V^{\ast\ast}\cong V: f\rightarrow \hat{f}$  which is continuous (see Proposition 32.22 in \cite{Zei90}). 
Hence for each $t\in [0,1]$, there exists $\psi_\epsilon(t) \in V$ such that 
\beq
J'_{\gamma_\epsilon(t)}[\psi_\epsilon(t)] = \|J'_{\gamma_\epsilon(t)}\|_{\ast}^2, \ \text{ and }\ \|\psi_\epsilon(t)\| = \|J'_{\gamma_\epsilon(t)}\|_{\ast}.
\eeq 
Since the map $t \to J'_{\gamma_\epsilon(t)}$ and the duality map are both continuous, we have the map $t \to \psi_\epsilon(t)$ also 
continuous.

For $h>0$, we consider the path
\be\label{rh}
\gamma_h(t) = \gamma_\epsilon(t) -h\rho_\epsilon(t)\psi_\epsilon(t) \in \Pat,
\ene
and let $t_h \in [0,1]$ be such that
\be\label{th}
J(\gamma_h(t_h)) = \max\limits_{t \in [0,1]} J(\gamma_h(t)) = F(\gamma_h).
\ene
On the one hand, by assumption \eqref{eke0b} we find:
\be\label{F1}
F(\gamma_h) \ge F(\gamma_\epsilon) - \sqrt{\epsilon}d(\gamma_\epsilon,\gamma_h) \ge J(\gamma_\epsilon(t_h)) - \sqrt{\epsilon}h\|\psi_\epsilon(t_h)\|.
\ene
On the other hand, we have 
\begin{eqnarray}\label{F2}
F(\gamma_h) &=& J\big(\gamma_\epsilon(t_h) -h\rho_\epsilon(t_h)\psi_\epsilon(t_h)\big) \notag \\
&=& J(\gamma_\epsilon(t_h)) - h\rho_\epsilon(t_h)J'_{\gamma_\epsilon(t_h)}\Big[\psi_\epsilon(t_h)\Big] + o(h).
\end{eqnarray}
Therefore from \eqref{F1} and \eqref{F2}, we find 
\beq
\rho_\epsilon(t_h)J'_{\gamma_\epsilon(t_h)}\Big[\psi_\epsilon(t_h)\Big]\le \sqrt{\epsilon}\|\psi_\epsilon(t_h)\| +o(1), \text{ as }\ h\to 0^{+}.
\eeq
This allows us to pass to the limit along a sequence $h_n \to 0^{+}$ with $t_{h_n} \to t_0\in [0,1]$. We claim that $t_0\in T_\epsilon$. Indeed, this 
follows from a general fact that if a sequence of continuous functions $f_n: [0,1]\to\R$ converges uniformly to a function $f$ 
(as $n\to\infty$), and if $f_n(t_n) =\max\limits_{t \in [0,1]}f_n(t)$, then along a subsequence, $t_n \to t_0$ with $f(t_0) =\max\limits_{t \in [0,1]}f(t)$.

And we conclude $J'_{\gamma_\epsilon(t_0)}\Big[\psi_\epsilon(t_0)\Big]\le \sqrt{\epsilon}\|\psi_\epsilon(t_0)\|$, and this gives
\beq
 \|J'_{\gamma_\epsilon(t_0)}\|_{\ast}^2 \le \sqrt{\epsilon} \|J'_{\gamma_\epsilon(t_0)}\|_{\ast}.
\eeq
Now the proof is complete.
\ep
%%%%%%%%%%%%%%%%%%%%
\subsection{Uniqueness}
In this subsection, we prove our main result:
\bt\label{add}
The Donaldson functional $\Dcal$ admits a unique {\cp} corresponding to its global minimum.
\et
\bp 
We assume by contradiction there are two distinct {\cp}s $P_1=(u_1,\eta_1)$ and $P_2=(u_2,\eta_2)$ for $\Dcal$ %, and consider the space of paths 
in $\V$, as in \eqref{path0}, we consider,
\be\label{path1}
\Pat = \big\{\gamma \in C^0([0,1], \V): \gamma(0) = P_1,\text{ and } \gamma(1) =P_2 \big\},
\ene
which defines a nonempty complete metric space, with
\beq 
d(\gamma_1,\gamma_2) = \max\limits_{t \in [0,1]}\|\gamma_1(t) - \gamma_2(t)\|_{\V}.
\eeq 
We will use Theorem ~\ref{Eke} and Lemma ~\ref{grad-D} on $(\V, \|\cdot\|_{\V})$, which is a uniformly convex Banach space (\cite{Cla36}). We take 
the functional $J = \Dcal$, and $F(\gamma) = \max\limits_{t \in [0,1]}\Dcal(\gamma(t))$. Set 
\be\label{c0}
c_0 = \max\big\{\Dcal(u_1,\eta_1), \Dcal(u_2,\eta_2) \big\},
\ene
so that, 
\beq
\max\limits_{t \in [0,1]}\Dcal(\gamma(t)) \ge c_0, \ \forall \gamma\in\Pat.
\eeq
Therefore, it is well defined:
\be\label{c}
c = \inf\limits_{\gamma\in\Pat}\max\limits_{t \in [0,1]}\Dcal(\gamma(t)) \ge c_0.
\ene
Note that, if the path $\gamma(t) = (u(t),\eta(t))$ lies in the space $\Pat$, then so does the path $\tilde{\gamma}(t) = (u(t),\eta(u(t)))$, where 
$\eta(u)\in W^{1,p}(X,E)$ is the map defined by Theorem ~\ref{du-map}. Indeed, since the map $\eta(u)$ is continuous, we see that $\tilde{\gamma}$ is 
also continuous, and we also check that $\tilde{\gamma}(0) = P_1$ and $\tilde{\gamma}(1) = P_2$, since by \eqref{3.23} we have $\eta(u(0)) =\eta(u_1) = \eta_1$ and 
$\eta(u(1)) =\eta(u_2) = \eta_2$. In addition, from Theorem ~\ref{du-map} we have:
\be\label{2.3}
\Dcal(\gamma(t)) \ge \Dcal(\tilde{\gamma}(t)), \ \forall t\in [0,1].  
\ene
We also emphasize that by applying Theorem ~\ref{Lem3} to the {\cp}s $P_1$ and $P_2$, there exist $\delta_0 > 0$ and $\epsilon_0 > 0$, such that 
\be\label{fence}
\Dcal(u,\eta(u)) \ge \Dcal(P_j) + \epsilon_0, \ \ \forall\ (u, \eta(u))\in\partial B_{\frac{\delta_0}{2}}(P_j), \ \, j=1,2;
\ene
and without loss of generality, we can assume further that $0 < \delta_0 \le \|P_1-P_2\|_{\V}$. Therefore for every $\gamma(t) = (u(t),\eta(t)) \in \Pat$, 
by continuity, we find $t_1, \ t_2\in [0,1]$ such that for $\tilde\gamma(t) = (u(t),\eta(u(t))) \in \Pat$, there holds:
\beq
\|\tilde{\gamma}(t_j)-P_j\|_{\V} = \frac{\delta_0}{2}, \ \ j =1, 2.
\eeq
So by \eqref{fence}, we have 
\be\label{fence2}
\max\limits_{t\in [0,1]}\Dcal(\gamma(t)) \ge \max\limits_{t\in [0,1]}\Dcal(\tilde{\gamma}(t)) \ge \Dcal(P_j) + \epsilon_0, \ \ j=1,2.
\ene
Therefore, for any $\gamma\in\Pat$, we have $\max\limits_{t\in [0,1]}\Dcal(\gamma(t)) \ge c_0 + \epsilon_0$ with suitable $\epsilon_0> 0$, and so the 
set 
\beq
T = \big\{\tilde{t} \in [0,1]: \max\limits_{t \in [0,1]}\Dcal(\gamma(t)) = \Dcal(\gamma(\tilde{t}))\big\}\subset\subset (0,1).
\eeq
Furthermore the value $c$ in \eqref{c} satisfies: $c \ge c_0+\epsilon_0$.

In view of \eqref{c} and \eqref{2.3}, for given $\epsilon > 0$, we find a continuous map $u_\epsilon^0: [0,1]\to H^1(X)$ such that 
\be\label{2.4a}
u_\epsilon^0(0) = u_1,\ \ \ \ u_\epsilon^0(1) = u_2,
\ene
and a path 
\be\label{2.4b}
\gamma_\epsilon^0(t) = (u_\epsilon^0(t), \eta(u_\epsilon^0(t))) \in \Pat
\ene
satisfying:
\be\label{2.5}
c \le\max\limits_{t \in [0,1]}\Dcal(\gamma_\epsilon^0(t)) < c+\epsilon.
\ene
At this point, we are in position to apply Theorem ~\ref{Eke} with the lower continuous map $F:\Pat\to[0,\infty)$, with 
$F(\gamma) = \max\limits_{t\in[0,1]}\Dcal(\gamma(t))$, {\wrt} the path $\gamma_\epsilon^0$ in \eqref{2.4b}. Thus, we obtain another path 
$\gamma_\epsilon(t)$ in $\Pat$, which in turn satisfies the assumptions of Lemma ~\ref{grad-D} with the functional $J = \Dcal$. Thus, we obtain 
$t_\epsilon \in (0,1)$, such that 
\be\label{tn}
c \le\Dcal(\gamma_\epsilon(t_\epsilon)) = \max\limits_{t \in [0,1]}\Dcal(\gamma_\epsilon(t)) < c+\epsilon,
\ene
and
\be\label{tn2}
d\big(\gamma_\epsilon(t_\epsilon),\gamma_\epsilon^0(t_\epsilon)\big) \le \sqrt{\epsilon},
\ene
and 
\beq
\|\Dcal'(\gamma_\epsilon(t_\epsilon))\|_{\ast} \le \sqrt{\epsilon}.
\eeq
Hence, as $\epsilon\to 0$, by \eqref{tn} and \eqref{tn2}, we have
\be\label{B}
\Dcal(\gamma_\epsilon(t_\epsilon)) \to c, \ \|\gamma_\epsilon(t_\epsilon)-\gamma^0_\epsilon(t_\epsilon)\|_{\V} \to 0,  \ \ \Dcal'(\gamma_\epsilon(t_\epsilon)) \to 0,
\ene
with $\gamma^0_\epsilon$ in \eqref{2.4b}.

Therefore, {\upto} $\epsilon_n \to 0$, applying Lemma ~\ref{apriori}, Proposition ~\ref{A}, and \eqref{B}, we can find sequences 
$(u_n,\eta_n)  := \gamma_{\epsilon_n}(t_{\epsilon_n})$ and $(\tilde{u}_n,\tilde{\eta}_n) :=\gamma_{\epsilon_n}^0(t_{\epsilon_n})$ such that 
\be\label{C}
u_n \xrightharpoonup{\ \ H^1} u, \ \text{ and } \ \tilde{u}_n \xrightharpoonup{\ \ H^1} u, \ \text{ as } n \to +\infty;
\ene 
also
\be\label{C2}
\tilde{\eta}_n = \eta(\tilde{u}_n) \xrightarrow{\ \ W^{1,p}(X,E)} \eta(u) =\eta,  \ {\eta}_n  \xrightarrow{\ \ W^{1,p}(X,E)} \eta, \ \text{ as } n \to +\infty.
\ene 
In other words, $(u_n,\eta_n)$ defines a {\PS} sequence for $\Dcal$, to which Lemma ~\ref{almostPS} applies. Hence, we conclude ({\upto}) that, 
$u_n\to u$ strongly in $H^1(X)$. In summarizing, we have established that, $(u_n,\eta_n) \to (u,\eta)$ strongly in $\V$ with $\Dcal(u,\eta) = c > c_0$ 
($c_0$ in \eqref{c0}), and $\Dcal'_{(u,\eta)} = 0$. Therefore $P_3 = (u,\eta)$ defines a {\cp} for $\Dcal$ different from $P_1$ and $P_2$.

We can apply Theorem ~\ref{Lem3} to $P_3$ to find suitable $\delta_0>0$, which we can always choose to satisfy: 
$0<\frac{\delta_0}{2}<\min\{\|P_j-P_3\|, j=1,2\}$, and $\epsilon_0 > 0$, such that 
\be\label{5.29}
\Dcal(u,\eta(u)) \ge \Dcal(P_3) + \epsilon_0, \ \ \forall\ (u, \eta(u))\in\partial B_{\frac{\delta_0}{2}}(P_3).
\ene
Since $\gamma^0_{\epsilon_n}(t_{\epsilon_n}) = (\tilde{u}_n,\tilde{\eta}_n) \to P_3$, as $n\to\infty$, while $\gamma^0_{\epsilon_n}(0)=P_1\not=P_3$ and 
$\gamma^0_{\epsilon_n}(1)=P_2\not=P_3$, by the continuity of $\gamma^0_{\epsilon_n}$, for $n$ sufficiently large, we find $\bar{t}_n \in (0,1)$ such that 
$\gamma^0_{\epsilon_n}(\bar{t}_n) \in \partial B_{\frac{\delta_0}{2}}(P_3)$, and so (by \eqref{5.29}) we have: 
$\Dcal(\gamma^0_{\epsilon_n}(\bar{t}_n)) \ge c+ \epsilon_0$. While in view of \eqref{2.5}, we also have:
\beq
\Dcal(\gamma^0_{\epsilon_n}(\bar{t}_n))\le \Dcal(\gamma^0_{\epsilon_n}(t_n)) \to c, \text{ as } n \to \infty,
\eeq
and we arrive at the desired contradiction.
\ep

%%%%%%%%%%%%%%%%
\section{Final Remarks}
%Holomorphic differentials on {\RS}s are central objects in {\Tt}, especially quadratic and cubic differentials, due to their close relationship with the study of 
%harmonic maps from closed surfaces. There are many recent important works also relating harmonic maps with representations of surface groups into 
%various character varieties, see for instance a recent survey \cite{Li19} and references within. In this section we explore some geometric 
%applications of our main results. 
%%%%%%%%%%%%%%%%%%%%%%%%%s
%\subsection{Minimal immersions of closed surfaces in {\htm}} 
As discussed in \cite{Uhl83} a {\mi} of $X$ into (a germ of) {\htm} $M\cong S\times\R$ (not necessarily complete) with a pullback metric $g = e^{2u}g_X$ 
and {\sff} $\text{II}_g = \Re(q)$ is governed by the Gauss-Codazzi equations:
\be\label{z2}	
   \left\{
   \begin{aligned}
     \frac{ \Delta u+1}{4}-\frac{\|q\|^2}{16}e^{-2u}-e^{2u}=0\ \ \ \ \ {\text{ on }} X,\\
      \db(q)=0,\ \ \ \ \ \ \ \ \ \ \ \ \  
   \end{aligned}
   \right.
\ene
expressing the Gauss consistency condition between intrinsic and extrinsic curvatures and the fact that $q$ defines a {\hqd} on $X$, namely $q \in \Cal_2(X)$.

More generally, if we concern with {\cmc} (CMC) immersions of {\mc} $c$, then the Gauss-Codazzi equations are modified accordingly as follows: 
\be\label{z3}	
   \left\{
   \begin{aligned}
     \frac{ \Delta u+1}{4}-\frac{\|q\|^2}{16}e^{-2u} -(1-c^2)e^{2u}=0\ \ \ \ \ {\text{ on }} X,\\
      \db(q)=0.\ \ \ \ \ \ \ \ \ \ \ \ 
   \end{aligned}
   \right.
\ene
If we attack \eqref{z3} by considering $q \in \Cal_2(X)$ fixed, it has been pointed out (\cite{HLT21}) that, by setting $\lambda= 1-c^2$, then the first 
equation in \eqref{z3} admits a solution if and only if $0 \le \lambda\le \lambda_0$ for a suitable $\lambda_0 = \lambda_0(q)$. Moreover we have multiple 
solutions when $0 < \lambda < \lambda_0$. Such a failure of one-to-one correspondence also prompted \cite{GU07} to take a different viewpoint inspired 
in some sense by the Higgs bundle approach pioneered by Hitchin (\cite{Hit87}). Hence, in \cite{GU07} the authors proposed to fix a class 
$[\beta] \in  \Hcal^{0,1}(X, E)$ with $E=T^{1,0}_X$ (i.e. $k =2$ in the notations above) and let $[\beta] = [\beta_0+\db\eta]$ where $\beta_0$ is the 
harmonic representative. For $\lambda= 1-c^2>0$ and $\eta$ a section of $T_X^{1,0}$ as above, now we aim to solve 
the system: 
\be\label{z4}	
   \left\{
   \begin{aligned}
     \frac{ \Delta u+1}{4}-\lambda e^{2u}-\|\beta_0+\db\eta\|^2e^{2u}=0\ \ \ \ \ {\text{ on }} X,\\
      \db(e^{2u}\ast_E(\beta_0+\db\eta))=0,\ \ \ \ \ \ \ \ \ \ \ \ \ 
   \end{aligned}
   \right.
\ene
and obtain a posteriori the {\hqd} $q = 4\ast_E(\beta_0+\db\eta)e^{2u}$, conveniently devised together with the metric $g = e^{2u}g_X$, and such that 
its real part identifies the traceless part of the {\sff} of the CMC immersion.

As before, a suitable change of variables reduces system \eqref{z4} to the Euler-Lagrange equations for the 
functional:
\be\label{df-lambda}
\Dcal_\lambda(u,\eta)=\int_X\{\frac14|\nabla u|^2- u+\lambda e^u+4\|\beta_0+\db\eta\|^2e^{u}\}dA.
\ene 
Obviously, for $\lambda> 0$, the functional $\Dcal_\lambda$ enjoys exactly the same properties of $\Dcal$ in \eqref{df}. So it admits a unique {\cp} 
$(u_\lambda,\eta_\lambda)$, and 
\beq
\Dcal_\lambda(u_\lambda,\eta_\lambda) = \min\limits_{\mathcal{W}}\Dcal_\lambda(u,\eta),  
\eeq
where $\mathcal{W}$ is defined in \eqref{domain} as the natural domain for the functional. As a consequence, we establish Corollary 3 as stated in the Introduction.
%%%%%%%%%%%%%%%%%%%%%%%%%%%%%%%%%%%%%

We conclude by pointing out that another advantage of our variational approach is that now we can aim to construct CMC immersions with 
$|c| =1$. Namely, we can try to see whether $(u_\lambda,\eta_\lambda)$ survives the passage to the limit as $\lambda \searrow 0$, in order 
to obtain a solution for problem \eqref{z4} when $\lambda = 0$. This is a nontrvial task, since for 
$\lambda = 0$, the functional $\Dcal_0 = \Dcal_{\lambda = 0}$ may not even be bounded from below. This occurs when $[\beta] = 0$ where we 
easily check that,  
\beq
\min\limits_{\mathcal{W}}\Dcal_\lambda = 4\pi(g-1)\ln(\lambda) \to -\infty, \ \ \text{as}\ \  \lambda \searrow 0,
\eeq
recall that $g$ is the genus of the surface. Hence a first interesting question is to understand the pair of data $(X,[\beta])$(if any) that yields to a functional 
$\Dcal_0$ bounded from below in $\mathcal{W}$. But, since $\Dcal_0$ is no longer coercive, in this case it is important to understand whether the 
minimum of $\Dcal_0$ is attained. To this purpose, it will be relevant to provide a detailed description about the asymptotic behavior of 
$(u_\lambda,\eta_\lambda)$ as $\lambda \searrow 0$, a task that becomes particularly delicate in case $\Dcal_0$ is unbounded from below in $\mathcal{W}$. 
To this purpose, one needs to elaborate a proper blow up analysis for sequence of solutions to \eqref{z4} as $\lambda \searrow 0$. This line of 
investigation will be expanded in future work.
%Nevertheless Theorem ~\ref{main} provides a parameter space $\mathcal{T}_g(S) \times \Hcal^{0,1}(X, E)$ which is of maximal dimension 
%$16g-16$, namely, as a corollary, we find:
%\begin{cor}
%Let $\mathcal{V}_0$ be the component of zero {\ti} in the representation space $Hom(\pi_1(S),SU(2,1))/SU(2,1)$, then
%there exists an open set of maximal dimension $16g-16$ (\cite{PP06, LM13}). Moreover, this open set can be parametrized
%by $(X, [\beta]) \in \mathcal{T}_g(S) \times \Hcal^{0,1}(X, E)$.\end{cor}

\bibliographystyle{amsalpha}
\bibliography{ref-para}
\end{document}